\documentclass[a4paper]{amsart}
\usepackage{amsthm,amsmath,amssymb,mathrsfs,url,amsaddr} 
\usepackage[numbers]{natbib}

\theoremstyle{definition}
\newtheorem{theorem}{Theorem}[section]
\newtheorem{lemma}[theorem]{Lemma}
\newtheorem{proposition}[theorem]{Proposition}
\newtheorem{corollary}[theorem]{Corollary}
\newtheorem{definition}[theorem]{Definition}
\newtheorem{remark}[theorem]{Remark}

\newcommand\DN{\newcommand}
\DN\lref[1]{Lemma~\ref{#1}}
\DN\tref[1]{Theorem~\ref{#1}}
\DN\pref[1]{Proposition~\ref{#1}}
\DN\sref[1]{Section~\ref{#1}}
\DN\dref[1]{Definition~\ref{#1}}
\DN\rref[1]{Remark~\ref{#1}} 
\numberwithin{equation}{section}

\newcounter{Const} 
\DN\Ct[1]{\refstepcounter{Const}c_{\theConst}\label{#1}}
\DN\Ctref[1]{c_{\ref{#1}}}

\DN\ts{\times}
\DN\ms{\medskip}
\DN{\ep}{\varepsilon}
\DN\st{\,;\,}
\DN{\limi}[1]{\lim_{#1\to\infty}} 	
\DN{\limz}[1]{\lim_{#1\to0}}

\DN\OD[2]{\frac{d #1 }{d #2}}
\DN\PD[2]{\frac{\partial #1 }{\partial #2}}
\DN\half{\frac{1}{2}}
\DN\map[3]{#1:#2 \to #3}

\DN\laweq{\stackrel{\mathrm{law}}{=}}

\DN\R{\mathbb{R}}
\DN\Rd{\mathbb{R}^d}
\DN\N{\mathbb{N}}
\DN\Q{\mathbb{Q}}
\DN{\Z}{\mathbb{Z}}
\DN\C{\mathbb{C}}

\DN\msaN{m _{\mathfrak{s},\alpha}^{N}}

\DN\Ss{\mathit{S}}
\DN\SSS{\mathsf{S}}
\DN\sss{\mathsf{s}}

\DN\muN{\mu ^{N}}
\DN\xx{\mathbf{x}}
\DN\XXX{\mathsf{X}}
\DN\WWW{\mathsf{W}}
\DN\III{\mathsf{I}}
\DN\yy{\mathbf{y}}
\DN\yyy{\mathsf{y}}
\DN\zz{\mathbf{z}}
\DN\uu{\mathbf{u}}
\DN\vv{\mathbf{v}}

\DN\E{\mathcal{E}}

\DN\ulab{\mathfrak{u}}
\DN\lab{\mathfrak{l}}

\DN\XX{\mathbf{X}}
\DN\YY{\mathbf{Y}}

\DN\Be{\mathrm{Be}}
\DN\mBe{\mathrm{mB}}

\DN\WN{W_{}^N}
\DN\WNp{W_{}^{N+1}}
\DN\Wn{W_{}^n}
\DN\intWN{\mathring{W}_{}^N}
\DN\Wi{W_{}^\infty}
\DN\WNNp{W_{}^{N,N+1}}
\DN\nWN{W_{\ge}^N}
\DN\nWNp{W_{\ge}^{N+1} }
\DN\nWn{W_{\ge}^n}
\DN\intnWN{\mathring{W}_{\ge}^N}
\DN\nWNN{W_{\ge}^{N,N}}
\DN\nWia{W_{\ge,\alpha}^\infty}

\DN\LaON{\Lambda_{\alpha, N}^\Omega}
\DN\LaONp{\Lambda_{\alpha, N+1}^\Omega}
\DN\LaNNp{\Lambda_{\alpha,N}^{N+1}}
\DN\LNNp{\Lambda_{N}^{N+1}}
\DN\LaNN{\Lambda_{\alpha,N}^{N}}
\DN\ELNNp{L_{N}^{N+1}}
\DN\ELiN{L_{N}^{\infty}}

\DN\MM{\mathcal{M}}
\DN\MMp{\mathcal{M}_{p}}
\DN\MMinv{\mathcal{M}_{p}^{\mathrm{inv}}}
\DN\MMerg{\mathcal{M}_{p}^{\mathrm{erg}}}
\DN\U{\mathbb{U}}

\DN\pp{\mathbf{p}}
\DN\diag{\mathrm{diag}}
\DN\rad{\mathfrak{rad}}
\DN\eval{\mathfrak{eval}}

\begin{document}
\title[Boundary processes associated with Laguerre and Pickrell diffusions]{Boundary Feller-Dynkin processes associated with Laguerre processes and Pickrell diffusions }
\author{ALEXANDER I. BUFETOV }
\address{Steklov Mathematical Institute of RAS, Moscow, Russia; Department of Mathematics and Computer Science, St. Petersburg State University,  St. Petersburg,  Russia; CNRS Institut de Math\'ematiques de Marseille, France}
\email{ bufetov@mi.ras.ru}

\author{YOSUKE KAWAMOTO}
\address{Graduate school of environmental, life, natural science and technology, Okayama University, Okayama, Japan}
\email{y-kawamoto@okayama-u.ac.jp}



\maketitle

\begin{center}
  
\it{
To Yakov Grigorievich Sinai on the occasion of his 90th birthday}

\end{center}
  
\begin{abstract}

In this paper, we construct a Feller-Dynkin boundary process by applying the method of intertwiners to the coherent family, introduced in our previous work, of Laguerre processes with a fixed parameter.
The corresponding boundary process is computed explicitly, and it turns out to be a deterministic dynamical system.
Since the Pickrell diffusions that leave the Pickrell measure invariant are coherent with respect to the same projective system, we also obtain the boundary process associated with the coherent family of Pickrell diffusions.
\end{abstract}

\section{Introduction}\label{s:1}

\subsection{The method of intertwiners}\label{s:11}
The aim of this paper is to give an explicit construction of boundary Feller-Dynkin processes associated with the projective system introduced in \cite{BuK25}, by applying the method of intertwiners.
Infinite-dimensional dynamics can be constructed in many examples \cite{BoO12, BoO13, Ols16,Cue18, Ass20}, and the boundary process is sometimes a stochastic infinite-dimensional diffusion \cite{Ols16b}.
Sometimes, however, a surprising phenomenon occurs, and in transition to infinite dimension the stochasticity disappears: the limit diffusion turns out to be a deterministic dynamical system \cite{Ass20}.
In this paper, we construct a projective system and verify that the resulting infinite-dimensional system is indeed deterministic, see \eqref{:54a} below.

We begin with a brief overview of this approach.
A detailed rigorous formulation is provided in \sref{s:2}.
A projective system $\{ E^N ,\ELNNp \}_{N\in\N } $ consists of a sequence of spaces $E^N $ and Markov kernels $\ELNNp $ from $E^{N+1} $ to $ E^{N}$.
The boundary of a projective system is a space $E^\infty$ equipped with a family of Markov kernels indexed by $N$ and satisfy  suitable consistency relations. 
In this sense, $E^\infty$ may be regarded as a projective limit of the spaces $\{E^N\}_{N\in\mathbb{N}}$.
Let $\{ T_t^N \}_{t\ge 0}$ be the Markov semigroup associated with a Markov process $\XX^N $ on $E^N $. 
If the intertwining relation 
\begin{align*}
 T_t^{N+1} \ELNNp =\ELNNp T_{t}^{N} \text{ for any }t\ge 0
\end{align*}
holds for each $N\in\N $, then we say that the family of Markov processes $\{ \XX^N \}_{N\in\N}$ is coherent (or consistent) with respect to $\{ E^N ,\ELNNp \}_{N\in\N } $.
Under certain additional conditions, a coherent family $\{\XX^N\}_{N\in\N}$ induces a corresponding Markov process $\XX^\infty $ on the boundary $E^\infty $.
Furthermore, if the kernels $\ELNNp  $ are Feller and the Markov processes $\{ \XX^N \}_{N\in\N}$ are Feller-Dynkin, then the boundary process $\XX^\infty $ is also Feller-Dynkin.

The method of intertwiners was introduced by Borodin and Olshanski \cite{BoO12}.
They constructed a Feller-Dynkin process on the boundary of the Gelfand-Tsetlin graph, which describes the branching of irreducible representations of the chain of unitary groups.
This approach was further applied to branching graphs associated with other groups by Borodin and Olshanski \cite{BoO13, Ols16} and by Cuenca \cite{Cue18}.
In their frameworks, the state spaces of the projective systems are discrete.
Assiotis \cite{Ass20} first applied the method of intertwiners to a continuous setting.
Let us explain the framework by Assiotis.

Let $\mathring{W}^N =\{ \xx = (x_i)_{i=1}^N \in\R^N \st x_1 < \cdots < x_N \}$ be the Weyl chamber and $\WN =\{ \xx \in\R^N \st x_1 \le \cdots \le x_N \}$ be its closure.
For $\xx=(x_i)_{i=1}^{N+1} \in \WNp $, we introduce the set
\begin{align*}
  \WNNp (\xx )& =\{\yy \in \WN \st x_1 \le y_1\le x_2\le \cdots \le y_N \le x_{N+1} \} 
.\end{align*}
For $\yy =(y_i)_{i=1}^N \in \WN $, we denote the Vandermonde determinant by
\begin{align*}
\Delta_{N}(\yy )=\prod_{1\le i <j\le N}(y_j- y_i)
.\end{align*}

In the particular case $\xx \in \mathring{W}^{N+1} $, we define the probability measure $\LNNp (\xx ,\cdot )$ on $\WN $ by the formula 
\begin{align}\label{:31g}
  \LNNp (\xx, d\yy ) = N! \cdot \frac{\Delta _N(\yy )}{\Delta _{N+1} (\xx _{})} \mathbf{1}_{ \WNNp (\xx )}(\yy) d\yy 
.\end{align} 
The kernel $\LNNp $ can be extended to a Feller kernel $ \WNp \dashrightarrow \WN $ \cite[Lemma 2.5]{Ass20}.
Thus, we have a projective system $\{\WN ,\LNNp \}_{N\in\N}$. 
The Dyson Brownian motions for the inverse temperature $\beta=2$, and its Ornstein-Uhlenbeck counterparts, form consistent families with respect to this projective system  \cite[Section 3.8]{AOW19}, \cite[Section 3]{War07} (see also \cite{GoS15, RaS18} for general $\beta $).
Furthermore, the family of diffusions leaving the Hua-Pickrell measures invariant is also coherent \cite{Ass20}.
Therefore, these coherent systems give rise to the associated boundary processes by a framework established in \cite{Ass20}. 

The Laguerre processes are also intertwined by $\LNNp $ \cite{Ass19a, AOW19}, but the parameter of the process varies with $N$; this is a shifted intertwining relation.
As a result, the Laguerre processes with a fixed parameter are not coherent families with respect to $\{\WN ,\LNNp \}_{N\in\N}$.
In \cite{BuK25}, we introduced a new kernel $\LaNNp $ which  intertwines the Laguerre processes with a fixed parameter.
We set $\nWN =\WN \cap [0,\infty )^N$.
Let $\alpha >-1 $.
For $\xx\in \mathring{W}_{\ge}^{N+1}:=\mathring{W}^{N+1} \cap (0,\infty )^{N+1}$, we define $\LaNNp (\xx, \cdot )$ as the probability measure on $\nWN $ given by
\begin{align}\notag &
\LaNNp (\xx, d\yy) 
\\& \label{:21b}
=N! (\alpha+1)_N  \frac{\Delta _N ( \yy )}{\Delta _{N+1} (\xx _{})} \prod_{k=1}^{N} \bigg( \mathbf{1}_{[x_{k-1} , x_{k+1}]}( y_k)  \int_{x_k \vee y_k }^{x_{k+1} \wedge y_{k+1} } \frac{y_k^{\alpha} }{z ^{\alpha+1}} d z  \bigg) d\yy 
.\end{align}
Here, we use the symbols $x_0=0$ and $y_{N+1} =\infty$ for notational convenience.
We also write $(x)_{n} =x(x+1)\cdots(x+n-1)$ for the shifted factorial.
Although \eqref{:21b} is defined  for $\xx \in\mathring{W}_{\ge}^{N+1} $, the kernel $\LaNNp $ extends to a Feller kernel $\nWNp \dashrightarrow \nWN $ \cite[Lemma~5, Proposition~8]{BuK25}.

In this paper, we apply the method of intertwiners to the projective system $\{ \nWN , \LaNNp \}_{N\in\N}$ for $\alpha \in \{0\} \cup \N $.
As a result, the boundary process corresponding to the Laguerre processes is obtained.
The boundary process turns out to be deterministic, and we find it explicitly; see \eqref{:54a} below.
Furthermore, we introduce diffusion processes leaving the Pickrell ensemble, the radial parts distribution of the Pickrell measure, invariant.
These dynamics are also intertwined by $\LaNNp $, and determine boundary processes. 

We note that the construction of the boundary and its associated processes relies on a matrix structure.
One of the key relations is the random matrix interpretation of the kernel $\LaNNp $ in terms of squared singular values.
This structure restricts $\alpha$ to the non-negative integers, and consequently the parameters $\alpha $ in constructing boundary processes are also restricted, although the corresponding intertwining relations are established for all $\alpha >-1$.
We expect the corresponding boundary result to remain valid for all $\alpha>-1$.

To state our main results, we give an informal description of the boundary of the projective system.
A precise definition will be given after the main results.
The boundary of the projective system $\{ \nWN , \LaNNp \}_{N\in\N}$ is given by
\begin{align}\label{:12d}
\Omega =\{(\boldsymbol{a }, \gamma )\st \boldsymbol{a }=(a_i)_{i\in\N}, a _1 \ge a _2 \ge \cdots \ge 0, \, \sum_{i\in \N} a _i \le \gamma  \}\subset \R^\N \ts \R
\end{align}
endowed with product topology.
By the definition of boundary, there exists a sequence of kernels $\LaON :\Omega \dashrightarrow \nWN $ such that, for each $N\in\N$, we have 
\begin{align*}
  \Lambda_{\alpha, N+1}^\Omega \LaNNp = \LaON
.\end{align*}

\subsection{The Laguerre Ornstein--Uhlenbeck process and the associated boundary process}

We first state our result on the Laguerre Ornstein--Uhlenbeck processes.
For $\alpha >-1$, we consider the $N$-dimensional stochastic differential equation
\begin{align} \label{:15a}
  dX_t^{N, i} &=\sqrt{2 X_t^{N,i} } dB_t^i +\Big( -X_t^{N, i} +\alpha +N +\sum_{j\neq i}^{N} \frac{ X_t^{N,i}+X_t^{N,j} }{X_t^{N,i} -X_t^{N,j} } \Big)dt
\end{align}
for $i=1,\ldots, N$.
Here, $(B^i )_{i=1}^N $ denotes the standard $N$-dimensional Brownian motion.
For $\alpha >-1$, the equation \eqref{:15a} has a unique strong solution in $\nWN $ for any starting point $\xx \in \nWN$ \cite[Theorem 2.2]{GrM14}.
We call the unique strong solution $\XX^N =(X^{N,i})_{i=1}^N$ the Laguerre Ornstein--Uhlenbeck process of parameter $\alpha $.
Let $\{ T_{\alpha, t}^N \}_{t\ge 0}$ be the Markov semigroup associated with $\XX^N$.

Note that equation \eqref{:15a} is the Ornstein--Uhlenbeck analogue of the non-colliding squared Bessel process, which have been extensively studied in the context of random matrix theory; see, e.g., \cite{KaT11}.
In fact, the non-colliding squared Bessel process describes the evolution of the squared singular values of a rectangular matrix-valued Brownian motion with complex-valued entries  \cite{Dem07,KoO01}; see also \cite{Bru91} for the real-valued case.

Define $C_{\infty }(\nWN )$ as the space of continuous functions on $\nWN $ vanishing at infinity.
We established the following intertwining relation.
\begin{proposition}\label{p:53}\cite[Theorem 1]{BuK25}
  Assume that $\alpha >-1$.
  Then, for any $N\in\mathbb{N} $, $f\in C_{\infty} (\nWN )$, and $t \ge 0$, we have the intertwining relation
  \begin{align}\label{:11b}
  T_{\alpha,t}^{N+1} \LaNNp f& =\LaNNp T_{\alpha,t}^{N}f
   .\end{align}
\end{proposition}

\begin{remark}
The intertwining relation \eqref{:11b} corresponds to the case of inverse temperature $\beta = 2$.
This relation was extended to general $\beta \ge 1 $ \cite{KaS26}.
\end{remark}

From \pref{p:53}, the family of Laguerre Ornstein--Uhlenbeck processes of parameter $\alpha $ is coherent with respect to $\{ \nWN , \LaNNp \}_{N\in\N}$.
Therefore, there exists a corresponding boundary process as follows.

\begin{theorem}\label{t:54}
There exists a unique Feller-Dynkin semigroup $\{ T_{t}^{\Omega } \}_{t\ge 0}$ on $\Omega $ such that, for any $\alpha \in \{0\} \cup \N$, $N\in\N$, $f\in C_{\infty} (W_{\ge}^N)$, and $t \ge 0$, we have
\begin{align*}
T_{t}^{\Omega } \LaON f=\LaON T_{\alpha ,t}^{N } f 
\end{align*}
and the boundary process is given by the formula \eqref{:54a}.
Let $\XX_t=(\boldsymbol{a } (t), \gamma (t))$ denote the process associated with $\{ T_{t}^{\Omega }  \}_{t\ge 0}$.
Then we have
\begin{align}\label{:54a}
     a_i (t)=a_i(0) e^{-t}, \quad \gamma(t)=1+ (\gamma (0)-1)e^{-t}
    .\end{align}
\end{theorem}
Note that, in general, a coherent family with respect to $\{ \nWN , \LaNNp \}_{N\in\N}$ yields boundary semigroups that may depend on $\alpha$. 
However, by Proposition~\ref{p:39}, the existence of an appropriate shifted intertwining relation implies that these semigroups are independent of $\alpha$. 
Since such a relation holds in the present setting, we denote 
the common boundary semigroup by $\{T_t^\Omega\}_{t\ge 0}$.

In \tref{t:54}, the boundary process is constructed for $\alpha\in \{0\}\cup\N $, although  \eqref{:11b} is valid for all $\alpha>-1$.
This is because the construction relies on a matrix structure as explained in \sref{s:14}.

The explicit formula \eqref{:54a} shows that $a_i(t)\to 0, \gamma(t)\to 1 $ as $t\to\infty$. 
Thus the boundary process converges to the point
$((0,0,\ldots),1)\in\Omega$, which can  be identified with the boundary measure associated with the consistent Laguerre ensembles. 
In this sense, the limit $t\to\infty$ may be interpreted as convergence to equilibrium on the boundary.

The boundary Laguerre Ornstein--Uhlenbeck process is deterministic and is given explicitly.
This is analogous to the results in \cite[Section 5.2]{Ass20}, where the boundary processes associated with the Dyson model and its Ornstein--Uhlenbeck counterpart are also deterministic and explicitly described.

\subsection{The Pickrell diffusion and the associated boundary process}
We identify another stochastic process with a rich integrable structure, called the Pickrell diffusion.
Actually, the family of Pickrell diffusions is also coherent with respect to the projective system $\{ \nWN , \LaNNp \}_{N\in\N}$.
As stated in \sref{s:15}, this diffusion is closely related to the Pickrell measures. 

For $\mathfrak{s} \in\R$ and $\alpha >-1$, we consider the $N$-dimensional stochastic differential equation
\begin{align}\notag 
dX_t^{N,i} &=\sqrt{2 X_t^{N,i}(1+X_t^{N,i} ) }dB_t^i 
\\& \label{:41a}\quad
+\Big( (2-2N-\mathfrak{s})X_t^{N,i} +\alpha+1 +\sum_{j\neq i}^{N}\frac{2X_t^{N,i} (1+X_t^{N,i})}{X_t^{N,i} -X_t^{N,j} } \Big)dt 
\end{align}
for $i=1,\ldots,  N$.
The well-posedness of this equation will be verified in \lref{l:32}.
Let $\XX ^N=(X^{N,i})_{i=1}^N $ be the unique strong solution to \eqref{:41a}, and  $\{ T_{\mathfrak{s},\alpha,t}^{N}  \}_{t\ge 0}$ be its semigroup.
We call $\XX ^N$ the Pickrell diffusion.
They are intertwined by $\LaNNp $ as follows.
\begin{theorem}\label{t:35}
Assume that $\mathfrak{s}\in \R $ and $\alpha >-1 $.
Then, for any $N\in\N$, $f\in C_{\infty} (W_{\ge}^N)$, and $t \ge 0$, we have the intertwining relation
\begin{align*}
T_{\mathfrak{s},\alpha ,t}^{N +1} \LaNNp  f=\LaNNp T_{\mathfrak{s},\alpha,t}^{N }f
.\end{align*}
\end{theorem}

We note that this intertwining relations are closely related to the $\beta$-corner Jacobi processes introduced in \cite{BoG15}; see \rref{r:410}.

\tref{t:35} shows that the family of Pickrell diffusions is coherent with respect to $\{ \nWN , \LaNNp \}_{N\in\N}$, which yields an associated boundary process.
Note that the boundary semigroup in \tref{t:12} is independent of $\alpha \in \{0\} \cup \N$ for the same reason as in \tref{t:54}.

\begin{theorem}\label{t:12}
 Assume that $\mathfrak{s}\in\R $.
  Then, there exists a unique Feller-Dynkin semigroup $\{ T_{\mathfrak{s},t}^{\Omega }  \}_{t\ge 0}$ on $\Omega $ such that, for any $\alpha \in \{0\} \cup \N$, $N\in\N$, $f\in C_{\infty} (W_{\ge}^N)$, and $t \ge 0$, we have
  \begin{align*}
  T_{\mathfrak{s} ,t}^{\Omega } \LaON f=\LaON T_{\mathfrak{s},\alpha ,t}^{N } f 
  .\end{align*}
Furthermore, if $\mathfrak{s}>-1$, the unique invariant probability measure for $\{ T_{\mathfrak{s},t}^{\Omega } \}_{t \ge 0 }$ is $m_\mathfrak{s} $, which will be defined in \sref{s:15}.
\end{theorem}

Unlike the boundary Laguerre process, the boundary Pickrell process constructed in \tref{t:12} is expected to be genuinely stochastic, rather than to degenerate into a deterministic process; see \cite[Section 10]{BoO13} and \cite[Example 6.2, Section 10]{Ols16b} for general arguments in related settings.
Furthermore, it is expected that the boundary Pickrell processes satisfies the infinite-dimensional stochastic differential equation
\begin{align} \label{:15f}
  dX_t^{i} &=\sqrt{2} X_t^{i} dB_t^i +\Big( -\mathfrak{s}X_t^{i} +\sum_{j\neq i}^{\infty} \frac{ X_t^{i} X_t^{j} }{X_t^{i} -X_t^{j} } \Big)dt
,\end{align}
for $i\in\N$.
Indeed, by using the identity
\begin{align*}
 \sum_{j\neq i}^{N} \frac{2x_i(1+x_i)}{x_i-x_j} = \sum _{j\neq i}^{N} \frac{2 x_i x_j +x_i +x_j }{x_i-x_j} +(N-1)(2x_i+1 )
,\end{align*}
the equation \eqref{:41a} can be written as
\begin{align} \notag
d X_t^{N,i} &=\sqrt{2X_t^{N,i} (1+X_t^{N,i}) } dB_t^i
\\&   \label{:41b}\quad
+\Big( -\mathfrak{s}X_t^{N,i} +N+\alpha +\sum_{j\neq i}^{N}\frac{2X_t^{N,i}X_t^{N,j} +X_t^{N,i} +X_t^{N,j} }{X_t^{N,i} -X_t^{N,j} } \Big)dt
.\end{align}
Thus, a formal scaling limit of \eqref{:41b} under the $N^{-2}$ scaling yields \eqref{:15f}. 
The scaling is motivated by the boundary embedding used in \sref{s:5}.
Note that equation \eqref{:15f} was derived in \cite{AsM} as the equation for a boundary process in a different setting, where both the underlying projective system and the approximating finite-particle systems differ from those considered in the present paper.
This suggests that the boundary process constructed in \cite{AsM} may be related to the boundary Pickrell process.
Moreover, the boundary process  in \cite{AsM} was shown to converge to equilibrium.
It is therefore natural to expect that the boundary Pickrell process has a similar long-time behaviour.

It would also be interesting to consider the Gibbs resampling property of the Laguerre and Pickrell diffusions, viewed as line ensembles, in the sprit of the Brownian Gibbs framework developed in \cite{CoH14}.
In our setting, Brownian bridges should be replaced by bridges of the corresponding one-dimensional diffusions.
Since \cite{CoH14} also treats the Gibbs property for an infinite line ensembles, it is natural to investigate analogous properties for the corresponding boundary processes.
For the boundary Laguerre process, the formula \eqref{:54a} suggests a deterministic version of Gibbs resampling, whilst it is less clear for the boundary Pickrell process.

\subsection{Ergodic decomposition for $\mathbb{U}(\infty) \ts \mathbb{U}(\infty) $-invariant probability measures}\label{s:14}

Let $M_{m,n}(\C)$ be the space of $m\times n$ matrices with complex entries, and for brevity write $M_{n}(\C):=M_{n,n}(\C)$.
We introduce the following subsets $H_{n}(\C), \mathbb{U}(n) \subset M_{n}(\C) $:  $H_{n}(\C) $ is the space of Hermitian matrices of size $n$, and $\mathbb{U}(n)$ is the space of unitary matrices of size $n$.
Let $M_{\N }(\C )$ be the projective limit $\varprojlim M_{N}(\C )$, the space of $\N \ts \N $ matrices.
Let $M_{\infty }(\C )$ denote the inductive limit  $\limi{N} M_{N}(\C )$, the space of $\N \ts \N $ matrices with finitely many non-zero entries.

Let $\mathbb{U}(\infty)= \limi{N}\mathbb{U}(N)$ denote the the inductive limit unitary group.
We say that a probability measure on $M_\N(\C)$ is $\U(\infty)\ts \U(\infty)$-invariant if it is invariant under the actions of multiplication by $\U(\infty )$ on the left and right.
We write $\mathcal{M}_p^{\text{inv} }$ for the set of all $\mathbb{U}(\infty)\ts \mathbb{U}(\infty) $-invariant probability measures on $M_\N(\C)$.
These invariant probability measures admit an ergodic decomposition, cf. \cite{Buf14}.
Let $\mathcal{M}_p^{\text{erg} } \subset \mathcal{M}_p^{\text{inv} } $ be the set of all ergodic $\mathbb{U}(\infty)\ts \mathbb{U}(\infty) $-invariant probability measures. 

For a probability measure $P$ on $M_{\N}(\C )$, we define the characteristic function of $P$ as the function on $M_{\infty }(\C )$ given by the formula
\begin{align*}
F_P(A)=\int_{M_\N(\C)}e^{i \mathfrak{Re}  \mathrm{Tr}(AX^*)}P (dX)
.\end{align*}
Note that, if $P\in \mathcal{M}_p^{\text{inv} }$, its characteristic function $F_P$ is uniquely determined by $F_P (A )$ for all diagonal matrices $A \in M_{\infty }(\C ) $.

The set $\mathcal{M}_p^{\text{erg} }$ can be characterised by $\Omega $, defined by \eqref{:12d}.
For $\omega=(\boldsymbol{a }, \gamma ) \in \Omega$, we write $\overline{\gamma}(\omega)=\gamma -\sum_{i\in\N}a_i$.
We define a function  $F_\omega $ on $\R$ by
\begin{align}\label{:12a}
F_\omega (\lambda )=e^{-4 \overline{\gamma}(\omega )\lambda^2 } \prod_{i\in\N } \frac{1}{1+4 a_i \lambda^2 }
.\end{align}
The next characterisation was first proved in \cite{Pic91} (see also \cite{Rab08}).

\begin{lemma}[\cite{Pic91}]\label{l:11}
There exists a bijective correspondence between $\Omega $ and $ \MMerg $.
We denote by $P_\omega \in \MMerg $ as the ergodic probability measure corresponding to $\omega\in\Omega$.
Furthermore, the characteristic function of $P_\omega \in \mathcal{M}_p^{\text{erg} }$ is given by, for $r_1,\ldots, r_n \in \R$,
\begin{align*}
\int_{M_\N(\C)}e^{i \mathfrak{Re} \mathrm{Tr}( \mathrm{diag} (r_1,\ldots,r_n,0,\ldots)X^*)} P_\omega (dX) =\prod_{j=1}^n F_\omega (r_j)
.\end{align*}
\end{lemma}

Based on \lref{l:11}, we will prove in \pref{p:23} that $\Omega$ is a boundary of $\{ \nWN , \LaNNp \}_{N\in\N}$.
Therefore, a coherence family of processes with respect to $\{ \nWN , \LaNNp \}_{N\in\N}$ leads to a boundary process on $\Omega$.
We next introduce a Markov kernel $\LaON  : \Omega \dashrightarrow \nWN $ as follows.

Define the map $\map{\mathfrak{rad}_{n}}{M_{m,n}(\C)}{\nWn }$ by the formula
\begin{align*}
\rad _{n} (X)=(\lambda_{1}(X), \ldots , \lambda_{n}(X))
,\end{align*}
where $(\lambda_{i}(X))_{i=1}^n$ denotes the eigenvalues of $X^{*}X$ arranged in non-decreasing order, which is referred as the radial part of $X$.
For a random matrix $X\in M_{m,n}(\C)$, let a probability measure $P_{\rad }^n[X]$ on $\nWn $ be the distribution of the radial part of $X$.

For $m_1\ge m_2, n_1\ge n_2 $, let $\pi_{m_2, n_2}^{m_1, n_1} : M_{m_1, n_1}  (\C) \to M_{m_2,n_2} (\C) $ be the natural projection sending an $m_1 \times n_1 $ matrix to its upper left $m_2 \times n_2$ corner. 
We employ the expression $\pi_{m_2,n_2}^{m_1}$ in place of $\pi_{m_2,n_2}^{m_1,n_1}$ if $m_1=n_1$, and use a similar symbol for $m_2=n_2$.
Furthermore, let $\map{\pi_{m,n }^\infty }{M_\N (\C) }{M_{m,n} (\C) }$ be the natural projection sending each $\N \ts \N$ matrix to its upper left $m \ts n$ corner.
We also use the symbol $\pi_{n}^{\infty }:=\pi_{n,n}^{\infty }$.

For $\alpha \in \{ 0\} \cup \N$, define the Markov kernel $\LaON : \Omega  \dashrightarrow \nWN $ by the formula
\begin{align}\label{:23b}
\LaON (\omega,\cdot)=(\mathfrak{rad}_{N} \circ \pi_{N+\alpha,N}^\infty )_{*} P_\omega
.\end{align}

\begin{remark}
Regarding the projective system $\{\WN , \LNNp \}_{N\in\N }$, the  explicit form of the boundary kernel to $\WN $ is given \cite[Equation (8)]{Ass20}.
It would be of interest to derive an explicit expression of $\LaON $.
\end{remark}

We note that intertwining relations for $\beta$-Laguerre Ornstein-Uhlenbeck processes via a $\beta$-analogues of $\LaNNp $ has been establish in \cite{KaS26}: these are $\beta$-analogues of the intertwining relation \eqref{:11b}.
It would therefore be of interest to construct the corresponding boundary $\beta$-Laguerre process.
At least for $\beta=1,4$, the above random matrix interpretation was extended to the $\beta$-analogues of $\LaNNp$ in \cite{KaS26}.
This suggests that an analogous construction of the boundary processes may also be possible in these case.

\subsection{The Pickrell measures}\label{s:15}

Let $dX$ be the Lebesgue measure on $M_{M,N} (\C)$. 
The Pickrell measure on $M_{M,N} (\C) $ is defined by 
\begin{align}\label{:13a}
P_{\mathrm{Pic},\mathfrak{s}}^{M,N} (dX ) \propto \det (1+X^{*} X)^{-M-N-\mathfrak{s}} \,dX
\end{align}
for a real number $\mathfrak{s}$.
We remark that this measure is finite if and only if $\mathfrak{s}>-1$.
Hereafter, when $\mathfrak{s}>-1$ holds, we regard $P_{\mathrm{Pic},\mathfrak{s}}^{M,N} $ as a probability measure on $M_{M,N} (\C) $  by normalising it with an appropriate constant.

For $\mathfrak{s} , \alpha >-1$, define the probability measure on $\nWN $ by the formula 
\begin{align*}
  \msaN (d \xx )=\frac{1}{Z_{\mathfrak{s},\alpha}^N} \Delta_N^2(\xx)  \prod _{k=1}^N x_k^\alpha  (1+x_k)^{-2N-\alpha-\mathfrak{s}}\, d\xx 
  ,\end{align*}
where $Z_{\mathfrak{s},\alpha}^N$ is the normalisation constant.
If $\alpha $ is a non-negative integer, the probability $\msaN $ is identical to the radial part of a random matrix according to the Pickrell measure $(\rad _N)_* P_{\mathrm{Pic},\mathfrak{s}}^{N+\alpha ,N}$.
Indeed, the push-forward measure of the Lebesgue measure $dX $ on $M_{N+\alpha, N}(\C)$ under the map $\rad_{N} $ is given by
\begin{align}\label{:14a}
(\rad_N)_* dX  = c_{N,\alpha}\Delta_{N}^2(\xx) \prod_{k=1}^{N} x_k^\alpha \, d\xx 
,\end{align}
for a positive constant $c_{N,\alpha}$ depending on $N$ and $\alpha$.
Thus, we call $ \msaN $ the Pickrell ensemble.
We remark that, for $\mathfrak{s}, \alpha >-1$, one can directly check that the Pickrell diffusion given by \eqref{:41a} is reversible with respect to the Pickrell ensemble $m_{\mathfrak{s},\alpha }^N $.
Hereafter in this subsection, we assume that $\mathfrak{s}>-1$.

The Pickrell measures have a consistency in the sense that
\begin{align*}
(\pi_{M_2,N_2}^{M_1,N_1})_* P_{\mathrm{Pic}, \mathfrak{s}}^{M_1,N_1} = P_{\mathrm{Pic}, \mathfrak{s}}^{M_2,N_2}
\end{align*}
from \cite[Proposition C.1]{Buf15}. 
Thereby, the Kolmogorov extension theorem yields a probability measure $P_{\mathrm{Pic}, \mathfrak{s}}$ on $M_\N(\C)$ such that 
\begin{align*}
(\pi_{M,N}^{\infty})_* P_{\mathrm{Pic}, \mathfrak{s}} = P_{\mathrm{Pic}, \mathfrak{s}}^{M,N}
.\end{align*}
We also call $P_{\mathrm{Pic}, \mathfrak{s}}  $ the Pickrell measure.

Clearly, $P_{\mathrm{Pic}, \mathfrak{s} }^{M,N}  $ is $\mathbb{U}(M)\ts \mathbb{U}(N) $-invariant under the actions by the unitary group of multiplication on the left and right.
Hence, $P_{\mathrm{Pic}, \mathfrak{s}} $ is also by the unitary group $\mathbb{U}(\infty) \ts \mathbb{U}(\infty) $.
Therefore, there exists a unique probability measure $m_{\mathfrak{s}} $ on $\Omega$ such that $P_{\mathrm{Pic}, \mathfrak{s}}  $ can be decomposed as 
\begin{align*}
P_{\mathrm{Pic}, \mathfrak{s}}  =\int_\Omega P_\omega d m_{\mathfrak{s}} (\omega) 
.\end{align*}
It is important to point out that $m_{\mathfrak{s}} $ is identified with the modified Bessel random point field $\mu_{\mBe ,\mathfrak{s}} $ in the following sense \cite[Propostion 1.10]{Buf15}.
We define the Bessel kernel, introduced in  \cite{TrW94}, given by the formula 
\begin{align*}
K_{\Be , \mathfrak{s}} (x, y)=\frac{\sqrt{x}  J_{\mathfrak{s}+1} (\sqrt{x} )  J_{\mathfrak{s}} (\sqrt{y} )-\sqrt{y}  J_{\mathfrak{s}+1} (\sqrt{y} ) J_{\mathfrak{s}} (\sqrt{x} )  }{2( x-y )}
\end{align*}
for $x,y>0$, where $ J_{\mathfrak{s}} $ denotes the Bessel function with parameter $\mathfrak{s}$.
The Bessel random point field $\mu_{\Be , \mathfrak{s}}$ is the determinantal random point field on $(0,\infty) $ associated with the Bessel kernel $K_{\Be , \mathfrak{s}}$ (more precisely, see \cite{Sos00}). 
The modified Bessel random point field $\mu_{\mBe ,\mathfrak{s}} $ is defined as the image measure of $\mu_{\Be , \mathfrak{s}}$ under the transformation $\sum_{i} \delta_{x_i} \mapsto \sum_{i} \delta_{4 /x_i } $.
Setting 
\begin{align*}
\Omega _{0}=\{(\boldsymbol{a }, \gamma )\in \Omega \st a_i > 0, \, \sum_{i\in \N}a_i = \gamma  \}
,\end{align*}
we have $m_{\mathfrak{s}} (\Omega_0)=1$.
Define the unlabelling map $\map{\ulab}{\Omega_{0} }{\SSS }$ by $\ulab((\boldsymbol{a },  \gamma) )=\sum_{i\in\N} \delta_{a_i}$, where $\SSS $ is the space of configurations on the open half-line $(0,\infty)$.
We then have
\begin{align*}
m_{ \mathfrak{s} } \circ \ulab^{-1}  =\mu_{\mBe , \mathfrak{s}}
.\end{align*}

\subsection{Organisation of this paper}
We introduce the general framework of the method of intertwiners in \sref{s:2}.
In \sref{s:3}, we apply the method of intertwiners to our projective system.
In \sref{s:4}, we prove the intertwining relations of Pickrell diffusions, and thus its boundary process is obtained.
In \sref{s:5}, we give an explicit description of the boundary Laguerre process.
In \sref{s:6}, we show that the kernel $\LaNNp $ arises as a scaling limit of a discrete kernel.

\section{Basic concepts of Feller boundaries and boundary processes}\label{s:2}
We first recall the method of intertwiners introduced in \cite{BoO12}. 
\subsection{Feller kernels and semigroups}\label{s:21}
We recall concepts of kernels and semigroups on Euclidean spaces; see \cite[Section 2]{BoO12} for a comprehensive discussion.
Let $E$ and $E'$ be Borel subsets of Euclidean spaces.
Consider a function $L(x,A)$ taking values in $[0,1]$, where $x\in E $ and $A$ is a Borel subset of $E'$.
We say that $L $ is a Markov kernel from $E$ to $E'$, denoted by $L : E \dashrightarrow E' $, if the following two conditions hold:
\begin{enumerate}
\item[$\bullet$]
 $L (x ,\cdot )$ is a probability measure on $E'$ for any $x\in E$.

\item[$\bullet$]
$L (\cdot, A) $ is a Borel function on $E$ for any Borel subset $A\subset E'$.
\end{enumerate}

Let $\mathcal{B}(E)  $ denote the Banach space of $\R $-valued bounded Borel functions with the sup-norm on $E$, and $C_\infty(E)$ denote the space of all continuous functions on $E$ vanishing at infinity.
Then, a Markov kernel $L : E \dashrightarrow E' $ defines a linear operator $\mathcal{B}(E') \to \mathcal{B}(E)$ by $(L f)(x)=\int_{E'}L (x, dy) f(y)$.
The kernel $L $ is called Feller if the induced operator satisfies $L (C_\infty(E')) \subset C_\infty(E) $.

Following \cite{BoO12}, we use the symbol $\mathcal{M}_p(E)$ to denote the space of all probability measures on $E$ equipped with the weak topology.
For $m \in \mathcal{M}_p(E)$, we write $m L (\cdot):=\int_{E} L (x,\cdot ) m (dx) \in \mathcal{M}_p(E')$.
Thus, a Markov kernel $L : E \dashrightarrow E' $ induces a map $\mathcal{M}_p(E) \to \mathcal{M}_p(E')$.

A Markov semigroup on $E$ is, by definition, a family of Markov kernels $\{ T_t\}_{t \ge 0}$ on $E$ such that $T_0 =1$ and $T_t T_s=T_{t+s}$ for all $t,s\ge 0$.
This semigroup is said to be Feller-Dynkin if the following two conditions hold:

\begin{enumerate}
\item[$\bullet$]
The induced map $T_t: \mathcal{B}(E) \to \mathcal{B}(E)$ satisfies $T_t ( C_\infty(E) ) \subset C_\infty(E) $ for any $t\ge 0$.

\item[$\bullet$]
The function $t\to T_t$ is strongly continuous in the sense that $\limz{t} T_t f =f $ in $C_\infty(E)$ holds for any $f\in C_\infty(E)$. 
\end{enumerate}

Note that Feller-Dynkin semigroups are sometimes simply called Feller semigroups.
To avoid confusion, we distinguish between these two in this paper. 
For us, the Feller property means just the property of mapping the space of continuous functions vanishing at infinity into itself; and we speak of the  Feller-Dynkin property if, additionally, the strong continuity holds.

\subsection{Feller boundary}
A projective system $\{E^N ,\ELNNp \}_{N\in\N}$ is, by definition, a sequence of Borel spaces $\{ E^N\}_{N\in \N}$ and Markov kernels $\ELNNp :E^{N+1} \dashrightarrow E^N$.
The family of induced maps $\map{\ELNNp }{\mathcal{M}_{p}(E^{N+1} )}{\mathcal{M}_{p}(E^N ) } $ forms a chain on the space $\prod_{N\in\N} \mathcal{M}_{p}(E^N ) $ with the product topology.
Define the set of all coherent probability measures by
\begin{align*}
\varprojlim \MMp (E^N )=\Big\{ \{ m^N \}_{N\in\N} \in \prod_{N\in\N} \MMp (E^N ) \st m^{N+1} \ELNNp = m^N \text{ for all }N \Big\}
.\end{align*}
\begin{definition}\label{d:22}
We say that a set $E^\infty $ is a boundary of a projective system $\{E^N, \ELNNp  \}_{N\in\N}$ if, for each $N\in\N$, there exists a Markov kernel $\ELiN : E^\infty \dashrightarrow E^N  $  such that the following conditions hold:
\begin{enumerate}
\item[(i)]
$L_{N+1}^\infty \ELNNp =\ELiN$.

\item [(ii)]
The induced map $\map{L^\infty}{\mathcal{M}_p(E^\infty ) }{\varprojlim \mathcal{M}_p( E^N )}$ defined by $L^\infty(P)=\{P \ELiN \}_{N\in\N}$ is a Borel isomorphism.
\end{enumerate}
Additionally, if $\{\ELNNp \}_{N\in\N}$ and $\{\ELiN \}_{N\in\N}$ are Feller kernels, then $E^\infty $ is called a Feller boundary.
\end{definition}

A boundary exists and is unique up to a Borel isomorphism \cite[Theorem 4.1.3]{Win85}.
This boundary coincides with the set of extremal coherent probability measures.
From the proof of \cite[Theorem 4.1.3]{Win85}, we have
\begin{align}\label{:21a}
E^\infty =\mathrm{Ex} \big(\varprojlim \mathcal{M}_{p} (E^N ) \big)
,\end{align}
where $\mathrm{Ex} $ denotes the set of extreme points.
By construction, for $\{m^N \}_{N\in\N}\in E^\infty$ and a Borel subset $A$ of $E^N$, the boundary kernel $L_K^\infty $ is given by 
\begin{align*}
L_K^\infty  \big(\{m^N \}_{N\in\N}, A  \big) =m^K (A )
.\end{align*}

\subsection{Intertwining relations and boundary processes}
A coherent family of processes produces a boundary process as follows. 
\begin{lemma}\label{l:26} \cite[Proposition 2.4]{BoO12}
For each $N\in\N$, let $T^N:= \{ T_t^N \}_{t \ge 0}$ be a Markov semigroup on $E_N$. 
Assume that the family $\{T^N\}_{N \in\N} $ is coherent with respect to the projective system $\{E^N, L_{N}^{N+1 }\}_{N\in\N }$ in the following sense: for $N\in \N$ and $t\ge 0$, we have
\begin{align}\label{:22c}
T_t^{N+1} \ELNNp =\ELNNp T_t^{N} 
.\end{align}
Then, there exists a unique Markov semigroup $T^{\infty} := \{ T_t^{\infty} \}_{t\ge 0 } $ on $E^\infty $ such that, for $t\ge 0$ and $N\in\N$,
\begin{align*}
T_t^\infty \ELiN =\ELiN T_t^{N}
.\end{align*}
Furthermore, if $E^\infty $ is a Feller boundary and $ T^N $ is a Feller-Dynkin semigroup for all $N\in\N$, then the semigroup $T^\infty$ is also Feller-Dynkin.
\end{lemma}

By \dref{d:22}, for any coherent family $\{ m^N \}_{N\in\N } \in \varprojlim \MMp (E^N ) $, there exists a unique element $m^\infty \in \MMp (E^\infty )$ such that 
\begin{align*}
m^\infty \ELiN =m^N 
.\end{align*}

\begin{lemma}\label{l:24} \cite[Section 2.8]{BoO12}
Suppose that the family of Markov semigroups $\{T^N\}_{N \in\N} $ is coherent with respect to $\{E^N, L_{N}^{N+1 }\}_{N\in\N }$.
Let $m^N $ be an invariant probability measure of $ T^N $. 
If $\{m^N\}_{N\in\N } \in \varprojlim \MMp (E^N ) $, the corresponding measure $m^\infty \in \MMp (E^\infty )$ is invariant for $T^\infty$.
Moreover, if $m^N$ is a unique invariant probability of $T^N $, then so is  $m^\infty$ for $T^\infty  $.
\end{lemma}

\begin{remark}
In the second statement in \lref{l:24}, we do not have to assume that $\{m^N \}_{N\in\N }$ is coherent.
In fact, the coherence follows from the assumption.
Actually, we have $m^{N+1} L_{N}^{N+1} T_t^N = m^{N+1} L_{N}^{N+1} $ from \eqref{:22c}.
Therefore, we obtain $m^{N+1} L_{N}^{N+1} = m^{N}$ from the uniqueness.
\end{remark}

\section{Boundary processes with respect to $\{ \nWN , \LaNNp \}_{N\in\N}$  }\label{s:3}
\subsection{Properties of the kernel $\LaNNp $}
Here we collect some facts about our kernel $\LaNNp $.
In this subsection, we assume $\alpha >-1$. 
We first observe that $\LaNNp $ given by \eqref{:21b} is decomposed into the kernel $\LNNp $ defined in \eqref{:31g} and a new kernel $\LaNN : \nWN \dashrightarrow \nWN $ defined as follows:
for $\xx \in \intnWN $, the kernel $\LaNN (\xx ,d\yy )$ is given by the formula 
\begin{align}\label{:32c}
  \LaNN (\xx , d\yy) = (\alpha+1)_{N} \bigg( \prod_{k=1}^{N}\frac{ y_k^\alpha }{ x_k^{\alpha+1}  }\bigg) \frac{\Delta_N (\yy)}{\Delta_N (\xx )}  \mathbf{1}_{\nWNN (\xx )}(\yy) d\yy
,\end{align}
where, for $\xx =(x_i)_{i=1}^{N} \in \nWN  $, 
\begin{align*}
  \nWNN (\xx ) :=\{\yy \in \nWN \st 0\le  y_1 \le x_1\le y_2\le \cdots \le y_N \le x_{N} \} 
.\end{align*}
Since $\LaNN (\xx, \cdot )$ is a probability measure on $\nWN $ \cite[Lemma 4]{BuK25}, the kernel $\LaNN  $ is Markov.
Furthermore, $\LaNN (\xx, \cdot )$ can be extended to all $\xx \in \nWN$ by \cite[Lemma 5]{BuK25}, and we obtain a kernel  $\LaNN : \nWN \dashrightarrow \nWN $.

\begin{lemma}\label{l:35} \cite[Proposition 7]{BuK25}
For any $\alpha >-1$ and $N\in\N$, we have
  \begin{align*}
  \LaNNp =\LNNp \LaNN 
  .\end{align*}
\end{lemma}

This decomposition is useful for establishing the intertwining relation in \tref{t:35}.
We next quote the Feller property of $\LNNp$, $\LaNN$, and $\LaNNp $.
\begin{lemma}\label{l:31}
  Suppose $\alpha >-1$ and $N\in\N$.
  Then, we have the following:
  \begin{align} \label{:31c}
  \LNNp f \in C_{\infty }(\nWNp ) \quad  &\text{ for any } f\in C_{\infty }( \nWN ),
  \\ \label{:31d}
  \LaNN f \in C_{\infty }(\nWN ) \quad &\text{ for any } f\in C_{\infty }( \nWN ),
  \\ \label{:31e}
  \LaNNp f \in C_{\infty }(\nWNp ) \quad &\text{ for any } f\in C_{\infty }( \nWN )
  .\end{align}
\begin{proof}
The property \eqref{:31c} was proved in \cite[Lemma 2.5]{Ass20}, and \eqref{:31d} and \eqref{:31e} were done in \cite[Proposition 8]{BuK25}.
\end{proof}
\end{lemma}

\subsection{The Feller Property of the kernel $\LaON$}
In the rest of this section, we assume that $\alpha \in\{0\}\cup \N$.
In this particular case, the kernels $\LaNN$ and $\LaNNp $ have interpretations in the context of random matrix theory.
We say that a random matrix $X_{m,n}\in M_{m,n}(\C)$ is $\mathbb{U}(m) \ts \mathbb{U}(n) $-invariant if $X_{m,n}\laweq V_{m} X_{m,n} U_{n}$ holds for any fixed  matrices $U_{n}\in \mathbb{U}(n)$, $V_{m}\in \mathbb{U}(m)$.

Let $V_{N+\alpha +1} \in \mathbb{U}(N+\alpha+1)$ be a Haar distributed random matrix.
Then, for any $\zz= (z_i)_{i=1}^{N} \in \nWN $, the kernel $\LaNN $ has the representation
\begin{align*}
\LaNN (\zz ,\cdot)= P_{\rad}^{N} [ \pi_{N+\alpha, N}^{N+\alpha+1}( V_{N+\alpha+1} ) \diag(\sqrt{z_1},\ldots, \sqrt{z_{N}}) ] 
\end{align*}
from \cite[Theorem 2.1]{KKS16}, where $P_{\rad}^{N}$ denotes radial parts distribution defined in \sref{s:14}.
Hence, for any $\U (N+\alpha+1)\ts \U(N)$-invariant random matrix $X_{N+\alpha+1,N}$, we have
\begin{align}\label{:31h}
 P_{\mathfrak{rad}}^{N} [ X_{N+\alpha+1,N} ] \Lambda_{\alpha, N}^{N}=P_{\mathfrak{rad}}^{N} [\pi_{N+\alpha,N}^{N+\alpha+1,N }( X_{N+\alpha+1, N})] 
.\end{align}
For this argument, see the proof of \cite[Theorem 2]{BuK25}:

It was also established in \cite[Theorem 2]{BuK25} that, for any $\mathbb{U}(N+\alpha+1) \ts \mathbb{U}( N+1)$-invariant random matrix $X_{N+\alpha+1,N+1} $, we have
\begin{align}\label{:31b}
P_{\rad}^{N+1} [ X_{N+\alpha+1,N+1} ] \LaNNp =P_{\rad}^{N} [\pi_{N+\alpha,N}^{N+\alpha+1, N+1}( X_{N+\alpha+1,N+1})]
.\end{align}
Furthermore, we obtain 
\begin{align*}
  \LaNNp(\xx , \cdot) =P_{\rad}^{N} \big[ \pi _{N+ \alpha ,N}^{N+\alpha+1, N+1}(V_{N+\alpha+1}  D_{N+\alpha+1, N+1}  U_{N+1}) \big]
\end{align*}
for any $\xx =(x_i)_{i=1}^{N+1}\in\nWNp$ \cite[Corollary 3]{BuK25}.
Here, $U_{N+1}  \in \mathbb{U}(N+1) $ and $V_{N+\alpha+1}  \in \mathbb{U}(N+\alpha+1) $ are Haar distributed independent random matrices, and $D_{N+\alpha+1, N+1} \in M_{N+\alpha+1, N+1}(\C)$ is a deterministic matrix given by
\begin{align*}
D_{N+\alpha+1,N+1}=
\begin{bmatrix}
\diag(\sqrt{x_1},\ldots,\sqrt{x_{N+1}} )  \\
\mathbf{0}_{\alpha \ts (N+1)}
\end{bmatrix}
.\end{align*}

Recall that the boundary kernel $\LaON $ is defined via matrix truncation \eqref{:23b}.
The above representations of $\LaNN$ and $\LaNNp $ (cf. \eqref{:32c}) now yield the following relationship between these kernels.

\begin{lemma}
For any $\alpha \in \{0\} \cup \N$ and $N\in\N$, we have the identities
  \begin{align}
\label{:24b}
  \Lambda_{\alpha+1,N}^\Omega \LaNN =\LaON 
,  \\  \label{:24a}
  \LaONp  \LaNNp =\LaON
  .\end{align}
   \begin{proof}
  The relations \eqref{:24b} and  \eqref{:24a} follow from \eqref{:23b} with \eqref{:31h} and \eqref{:31b}, respectively.

\end{proof}
\end{lemma}

The relation \eqref{:24a} is the consistency relation required in condition (i) of \dref{d:22}.
The relation \eqref{:24b} is needed in \pref{p:39} to show that the boundary processes constructed for different for values of $\alpha $ coincide.

Let $\omega_n,\omega \in\Omega $ be such that $\limi{n}\omega_n=\omega$.
Then, by the definition of characteristic functions \eqref{:12a}, we have 
$$\limi{n}\prod_{i=1}^N F_{\omega_n}(x_i )=\prod_{i=1}^N F_{\omega}(x_i )$$
for any $(x_1,\ldots,x_N)\in\nWN $.
Therefore, from the L\'evy continuity theorem with \eqref{:23b}, we have the continuity of $\LaON f$: for any $f\in C_b(\nWN )$,
\begin{align}\label{:33c}
\limi{n } \LaON f(\omega_n)    =\LaON f(\omega)
.\end{align}

To prove the Feller property of $\LaON $, following the strategy of \cite[Lemma 2.12]{Ass20}, we derive an analytic representation of the kernel.
\begin{lemma}\label{l:28}
Suppose that $\alpha \in \{0\} \cup \N$.
Let $f\in \mathcal{S}(\nWN )$ be a Schwartz function on $\nWN $.
Then, for any $\omega\in\Omega $, we obtain 
\begin{align}\label{:34a}
\LaON f (\omega)=c_{N,\alpha} \int _{\nWN }  \Delta_{N}^2(\xx) \Big(\prod_{k=1}^{N}x_k^\alpha \Big) F_{\omega}(x_1)\cdots F_{\omega}(x_N) \overline{\widehat{f}(\xx) } \, d\xx  
,\end{align}
where $c_{N, \alpha }$ is the constant same as in \eqref{:14a}.
\begin{proof}
We first consider the case where $\omega \in\Omega $ satisfies $\overline{\gamma }(\omega )>0$.
In this case, the Gaussian factor in the characteristic function \eqref{:12a} shows that the finite-dimensional projection $(\pi_{N+\alpha, N}^\infty )_* P_\omega $ is the convolution of a
probability measure with a non-degenerate Gaussian measure.
Hence, this projection has a probability density function denoted by $P_\omega^{N+\alpha, N} (X)$ for the Lebesgue measure $dX$ on $M_{N+\alpha ,N}(\C )$.
Define a $\U (N+\alpha) \ts \U(N)$-invariant function  $G_\omega $ on $M_{N+\alpha,N}(\C)$ such that $G_\omega (X)=F_\omega(x_1)\cdots F_\omega(x_N )$, where $\rad_N(X)=(x_1,\ldots,x_N)$.
Then, the Fourier transform of $P_\omega^{N+\alpha, N} (X)$ is given by $G_\omega (X)$.
Therefore, by the Plancherel theorem, we obtain 
\begin{align*}
  \LaON f (\omega)&=\int_{M_{N+\alpha, N}(\C)} f(\rad_N (X ))  P_\omega^{N+\alpha, N} (X) \, dX
  \\&
  =\int_{M_{N+\alpha, N}(\C)} G_\omega (X) \overline{\widehat{f}(\rad_N (X ))} \, dX
.\end{align*}
This is identical to the right-hand side of \eqref{:34a} from \eqref{:14a}.
Thus, we have shown \eqref{:34a} if $\overline{\gamma }(\omega ) >0$.

For general $\omega \in \Omega $, choose a sequence $\{\omega_{n}\}_{n\in\N} \subset \Omega $ such that $\limi{n}\omega_{n}=\omega$ and $\overline{\gamma}(\omega_n)>0$ for all $n$.
Then \eqref{:34a} holds for each $\omega_n$.
Letting $n\to \infty $ in this identity, the left-hand side converges to $\LaON f(\omega )$ by \eqref{:33c}, whilst the right-hand side converges to the right-hand side of \eqref{:34a} by the dominated convergence theorem.
Thus, \eqref{:34a} holds for any $\omega \in \Omega $, which completes the proof. 
\end{proof}
\end{lemma}
  
\begin{lemma}\label{l:25}
For each $\alpha \in \{0\} \cup \N$ and $N\in\N$, the kernel $\LaON $ is Feller.
\begin{proof}
Let $f\in C_{\infty} (\nWN)$.
The continuity of $\LaON f$ on $\nWN $ has been proved in \eqref{:33c}.
It remains to prove that $\LaON f$ vanishes at infinity.
It suffices to verify this for $f \in \mathcal{S}(\nWN )$ since $\mathcal{S}(\nWN ) $ is dense in $C_\infty (\nWN )$.

Let $\ep >0$.
Since $\hat{f}\in \mathcal{S}(\nWN )$, there exists $R$ such that
\begin{align*}
  c_{N,\alpha} \int _{\nWN (R ) } \Big| \Delta_{N}^2(\xx) \Big(\prod_{k=1}^{N}x_k^\alpha \Big)  F_{\omega }(x_1)\cdots F_{\omega }(x_N) \overline{\widehat{f}(\xx) } \Big| \, d\xx < \ep
,\end{align*} 
where  $\nWN (R ):= \nWN \cap ( [0,R]^N)^c$.
Thus, \lref{l:28} yields
\begin{align}\notag
  | \LaON f (\omega ) | &<  \ep + c_{N,\alpha} \int _{\nWN \cap [0,R ]^N}\Big|  \Delta_{N}^2(\xx) \Big(\prod_{k=1}^{N}x_k^\alpha \Big)  F_{\omega }(x_1)\cdots F_{\omega }(x_N) \overline{\widehat{f}(\xx) } \Big|d\xx 
  \\& \label{:35a}
  \le \ep + c_{N,\alpha} \Big( \int_{0}^R |x^\alpha F_{\omega} (x) | dx \Big)^N \sup_{\xx \in \nWN \cap [0,R ]^N}  |\Delta_{N}^2(\xx)  \overline{\widehat{f}(\xx) } | 
.\end{align} 

We now suppose that $\omega \to \infty $, which is reduced to the following three cases by the definition of $\omega$: (i) $\overline{\gamma }(\omega) \to \infty $, (ii)  $a_1 \to \infty $, (iii)  $a_1 $ remains bounded, and $\sum_{i\in\N} a_i \to\infty $.
Because we observe that 
\begin{align*}
\prod_{i \ge 2 } \frac{1}{1+4 a_i x^2}
\le \prod_{i\ge 2 } \exp \Big(- \frac{4 a_i x^2}{1+4 a_i x^2} \Big)
\le  \exp \Big( -\frac{4 x^2}{1+4 a_1 x^2}\sum_{i \ge 2 } a_i \Big)
,\end{align*} 
we obtain 
\begin{align*}
F_\omega (x) \le \frac{e^{-4\overline{\gamma}(\omega ) x^2}}{1+4 a_1 x^2} \exp \Big( -\frac{4 x^2}{1+4 a_1 x^2}\sum_{i \ge 2 } a_i \Big)
.\end{align*} 
Therefore, combining this with \eqref{:35a}, we conclude that $\LaON f$ vanishes at infinity in any of the cases (i), (ii), and (iii). 
Thus, the proof is complete.
\end{proof}
\end{lemma}

\subsection{ Feller boundary of the projective system $\{\nWN ,\LaNNp \}_{N\in\N} $}
In this subsection, we prove that $\Omega $ is a Feller boundary of the system $\{\nWN ,\LaNNp \}_{N\in\N} $. 
As before, we assume that $\alpha \in \{0\}\cup \N $ in this subsection.
Recall that $\MMp (\nWN )$ denotes the set of all probability measures on $\nWN $.
Because the kernel $\LaNNp  $ induces a map $\MMp (\nWNp ) \to \MMp (\nWN ) $ as stated in \sref{s:21}, 
the family of kernels $\{ \LaNNp  \}_{N\in\N }$ gives the chain
\begin{align*}
\cdots  \to \MMp (W_{\ge}^{N+1} ) \to \MMp (W_{\ge}^N ) \to \cdots \to \MMp (W_{\ge}^2) \to \MMp (W_{\ge}^1)
.\end{align*}
Define the set of all coherent probability measures with respect to this chain by
\begin{align*}
\varprojlim \mathcal{M}_{p,\alpha } (\nWN )=\big\{ \{ m^N \}_{N\in\N} \in \prod_{N\in\N} \mathcal{M}_{p}(\nWN )  \st m^{N+1} \LaNNp = m^N \big\}
.\end{align*}
The projective system $\{\nWN ,\LaNNp \}_{N\in\N}$ has the boundary $\nWia $ given by
\begin{align*}
\nWia =\mathrm{Ex} \big(\varprojlim \mathcal{M}_{p,\alpha} (\nWN ) \big)
\end{align*}
from \eqref{:21a}. 
The boundaries seem to depend on $\alpha $ at this stage, but we shall show that these all are identical to $\Omega$ and independent of $\alpha $.

For any $P\in \mathcal{M}_{p}^{\mathrm{inv}}$, we have 
\begin{align*}
  \big( (\mathfrak{rad}_{N+1} \circ \pi_{N+\alpha+1, N+1}^\infty )_{*} P \big)\LaNNp =(\mathfrak{rad}_{N} \circ \pi_{N+\alpha, N}^\infty )_{*} P
\end{align*}
by \eqref{:31b}.
Hence the family $ \big\{(\mathfrak{rad}_{N} \circ \pi_{N+\alpha, N}^\infty )_{*} P \big\}_{N\in\N}$
belongs to $\varprojlim \mathcal{M}_{p,\alpha}(\nWN)$. 
Therefore, the following map is well-defined.
\begin{definition}
Define a map $\map{}{\mathcal{M}_{p}^{\mathrm{inv}}}{\varprojlim \mathcal{M}_{p, \alpha} ( \nWN )}$ by 
\begin{align*}
  \Phi_\alpha(P)=\{(\mathfrak{rad}_{N} \circ \pi_{N+\alpha, N}^\infty )_{*} P \}_{N\in\N }
.\end{align*}
\end{definition}

\begin{lemma}\label{l:43}
The map $\map{ \Phi_{\alpha} }{\mathcal{M}_{p}^{\mathrm{inv}}}{\varprojlim \mathcal{M}_{p, \alpha} ( \nWN )}$ is a bijection.
\begin{proof}
We first note that, a $\U(N+\alpha ) \ts \U(N)$-invariant probability measure on $M_{N+\alpha, N}(\C) $ is characterized by its radial parts distribution $m^N \in\mathcal{M} _p( \nWN )$ as follows.
For $m^N \in\mathcal{M} _p( \nWN )$, let $V_{N+\alpha}, U_{N} ,D_{N+\alpha,N} (m^N)$ be independent random matrices such that $V_{N+\alpha} \in \mathbb{U}(N+\alpha) , U_{N}  \in \mathbb{U}(N) $ are Haar distributed and 
\begin{align*}
D_{N+\alpha,N}(m^N) =
\begin{bmatrix}
	\mathrm{diag}(\sqrt{x_1},\ldots,\sqrt{x_{N}} )  \\
\mathbf{0}_{\alpha \times N}
\end{bmatrix}
\in M_{N+\alpha, N}(\C)
,\end{align*} 
where the distribution of $(x_1,\ldots, x_{N})$ is given by $m^N$.
We set $X_{N+\alpha, N} (m^N) := V_{N+\alpha} D_{N+\alpha,N} (m^N) U_{N}$.
By construction, $X_{N+\alpha, N} (m^N) $ leads to a $\U(N+\alpha ) \ts \U(N)$-invariant probability measure on $M_{N+\alpha, N}(\C) $.
Furthermore, by a similar reason as in \cite[Lemma 2.4]{Def10} (see also \cite[Equation (30)]{BuK25}), a random matrix which has radial parts law $m^N$ is identical to $ V_{N+\alpha} D_{N+\alpha,N} (m^N) U_{N}$ in law.
Thus, a $\U(N+\alpha ) \ts \U(N)$-invariant probability is determined by its radial parts distribution.

For $m^{N+1} \in \MMp (\nWNp )$, let $X_{N+\alpha +1, N+1} (m^{N+1}) $ be as above.
Clearly, the random matrix 
$$\pi_{N+\alpha , N}^{N+\alpha +1, N+1} (X_{N+\alpha +1, N+1} (m^{N+1}) )$$
is $\U(N+\alpha ) \ts \U(N)$-invariant.
Furthermore, by \eqref{:31b}, its radial parts has law $m^{N+1}\LaNNp $.
Therefore, we have 
\begin{align}\label{:43e}
\pi_{N+\alpha , N}^{N+\alpha +1, N+1} (X_{N+\alpha +1, N+1} (m^{N+1}) ) \laweq X_{N+\alpha, N} (m^{N+1}\LaNNp)
.\end{align}

For $\{ m ^N \}_{N\in\N } \in \varprojlim \mathcal{M}_{p,\alpha} ( \nWN ) $, let $P^{N+\alpha, N} $ be the distribution of the random matrix $X_{N+\alpha, N} (m^N )$.
Then, from \eqref{:43e} and the consistency $m^{N+1}\LaNNp =m^N$, the sequence $\{P^{N+\alpha, N} \}_{N \in\N}$ satisfies the consistency
\begin{align*}
(\pi_{N+\alpha, N}^{N+\alpha+1, N+1})_* P^{N+\alpha +1 ,N+1}= P^{N+\alpha, N}
.\end{align*}
Hence, by the Kolmogorov extension theorem, there exists a probability measure $P$ on $M_{\N }(\C) $ such that 
\begin{align*}
  (\pi_{N+\alpha, N}^{\infty})_* P= P^{N+\alpha, N}
.\end{align*}
Since $P^{N+\alpha,N}$ is $\U(N+\alpha)\ts \U(N)$-invariant, we have $P\in\MMinv $.
By construction, we obtain $\Phi_\alpha (P) =\{ m ^N \}_{N\in\N}$, which implies that $\Phi_\alpha$ is a surjection.

To prove injectivity, let $P,Q\in \mathcal{M}_{p}^{\mathrm{inv}}$ satisfy $\Phi_\alpha (P)=\Phi_\alpha (Q)$.
Then, for the same reason as above that the $\U(N+\alpha )\ts \U(N)$-invariant distribution on $M_{N+\alpha, N}(\C )$ is characterised by its radial parts law, we have $(\pi_{N+\alpha, N}^{\infty})_* P =(\pi_{N+\alpha, N}^{\infty})_* Q$.
Thus all finite-dimensional projections coincide, and hence $P=Q$.
Therefore $\Phi_\alpha$ is injective, and we complete the proof.
\end{proof}
\end{lemma}

By \eqref{:24a}, the map $\Lambda_\alpha^{\Omega }:\Omega \to \varprojlim \mathcal{M}_{p, \alpha} ( \nWN ) $ defined by $$  \Lambda_\alpha^{\Omega } (\omega )  =\{\LaON (\omega, d\xx )  \}_{N\in\N} $$ is well defined.
Define the induced map $\map{\tilde{\Lambda}_\alpha^{\Omega }}{\MMp (\Omega )}{  \varprojlim \mathcal{M}_{p, \alpha} ( \nWN )  } $ by
\begin{align*}
\tilde{\Lambda}_\alpha^{\Omega } (m)& =\Phi_\alpha \Big(\int P_\omega dm (\omega ) \Big) 
.\end{align*}  
The map $\tilde{\Lambda}_\alpha^{\Omega } $ is an affine bijection from \lref{l:11} and \lref{l:43}. 
We remark that $\Lambda_\alpha^\Omega $ gives a bijection between  $\Omega $ and $\mathrm{Ex}\big(\varprojlim \mathcal{M}_{p, \alpha} ( \nWN ) \big)$.

We now prove that $\Omega $ is a Feller boundary, following the arguments in  \cite[Section 2.2]{Ass20} and \cite[Theorem 3.1]{BoO12}.
\begin{proposition}\label{p:23}
For $\alpha \in \{0\}\cup \N$, the Feller boundary of $\{  \nWN ,\LaNNp \}_{N\in\N}$ is given by $\Omega $ with the family of kernels $\{\LaON \}_{N\in\N}$.
In particular, the boundary space $\Omega $ is independent of $\alpha$.
\begin{proof}
From \lref{l:25}, the map $\map{\Lambda_\alpha^{\Omega } }{\Omega }{ \mathrm{Ex}\big(\varprojlim \mathcal{M}_{p, \alpha} ( \nWN ) \big)}$ is continuous, which implies that $\tilde{\Lambda}_\alpha^{\Omega } $ is Borel.
To show that the inverse of $\tilde{\Lambda}_\alpha^{\Omega } $ is also Borel, we can use a result in \cite{Mac57}.
Actually, since $\Omega $ is a standard Borel space, so does $\MM_p(\Omega )$.
Therefore, from \cite[Theorem 3.2]{Mac57}, we conclude that  $\tilde{\Lambda}_\alpha^{\Omega } $ is a Borel isomorphism.
Combining this with \eqref{:24a}, we prove that $\Omega $ is a boundary.
The Feller assertion comes from \lref{l:31} and \lref{l:25}.
\end{proof}
\end{proposition}

The method of intertwiners yields boundary processes on $\Omega$ as follows.

\begin{proposition}\label{p:39}
Suppose $\alpha \in \{0\}\cup \N$. 
Assume that a family of Markov semigroups $\{T_{\alpha }^N\}_{N\in\N }$ on $\nWN $ satisfies, for any $N\in\N$ and $t\ge 0$, 
\begin{align}\label{:39a}
T_{\alpha, t}^{N+1} \LaNNp =\LaNNp T_{\alpha, t}^{N}
.\end{align}
Then, there exists a unique Markov semigroup $T_{\alpha  }^\Omega :=\{T_{ \alpha ,t }^\Omega \}_{t\ge 0}$ on $\Omega $ such that
\begin{align}\label{:39c}
  T_{\alpha, t}^{\Omega } \LaON =\LaON T_{\alpha, t}^{N}
.\end{align}
Furthermore, if $T_{\alpha }^N$ is Feller-Dynkin, so is $T_{\alpha }^\Omega $.
Additionally, if 
\begin{align}\label{:39b}  
 T_{\alpha+1, t}^{N} \LaNN =\LaNN T_{\alpha, t}^{N}
\end{align}
holds, then $T_{\alpha+1, t}^{\Omega }=T_{\alpha, t}^{\Omega }$.
\begin{proof}
The first assertion follows from \lref{l:26}.
Using \eqref{:39c} and \eqref{:39b}, we obtain $T_{\alpha+1, t}^{\Omega } \Lambda_{\alpha+1,N}^{\Omega } \LaNN  =\Lambda_{\alpha+1, N}^{\Omega}\LaNN T_{\alpha, t}^{N}$, which implies $T_{\alpha+1, t}^{\Omega } \LaON =\LaON T_{\alpha, t}^{N}$ from \eqref{:24b}. 
It follows that $T_{\alpha+1, t}^{\Omega }=T_{\alpha, t}^{\Omega }$ by uniqueness.
Thus, we complete the proof.
\end{proof}
\end{proposition}

\section{Construction of the boundary Feller-Dynkin processes associated with Pickrell diffusions}\label{s:4}

\subsection{Pickrell diffusions and the Karlin-MacGregor semigroups}\label{s:41}
For $\mathfrak{s}\in\R $ and $\alpha >-1$, we consider a diffusion on $[0,\infty) $ associated with the generator
\begin{align*}
L_{\mathfrak{s},\alpha ,x}^{(N)} :=L_{\mathfrak{s},\alpha }^{(N)}:=x(1+x) \frac{d ^2}{d x^2}+\{ (2-2N-\mathfrak{s})x +(\alpha+1) \}\OD{}{x}
\end{align*}
with the following boundary conditions (see, for example, \cite{EtK86} for a detailed discussion of boundary conditions): the point $\infty$ is a natural boundary, and the  origin is an entrance boundary for $\alpha\ge 0$ and a regular boundary for $-1<\alpha<0$, in which case we impose the reflecting boundary condition.
The stochastic differential equation associated with $L_{\mathfrak{s},\alpha }^{(N)}$ is
\begin{align*}
dX_t= \sqrt{2X_t(1+X_t)} dB_t +\{(2-2N-\mathfrak{s})X_t +(\alpha+1)\} dt 
.\end{align*}

\begin{lemma}
The Vandermonde determinant $\Delta _N(\mathbf{x})$ is an eigenfunction of the second-order operator $\sum_{ i=1}^{N} L_{\mathfrak{s}, \alpha ,x_i}^{(N)} $ with eigenvalue $$\lambda _\mathfrak{s}^N =\frac{N(N-1)(-4N+2-3\mathfrak{s})}{6}.$$
\begin{proof}
We can check this lemma by a direct computation in \cite{Kon05}.
\end{proof}
\end{lemma}

Let $p_{\mathfrak{s},\alpha ,t}^{(N)}(x,y)$ be the transition density of the diffusion associated with $L_{\mathfrak{s},\alpha }^{(N)}$.
For $(t,\xx, \yy)\in (0,\infty)\ts \intnWN \ts \nWN$, consider the Karlin-McGregor transition density of $N$ particle $L_{\mathfrak{s},\alpha}^{(N)}$-diffusions h-transformed by $\Delta_N(\mathbf{x})$: 
\begin{align*}
\pp_{\mathfrak{s},\alpha ,t}^N(\mathbf{x},\mathbf{y})
:=e^{-\lambda_{\mathfrak{s}}^N t} \frac{\Delta_N(\yy ) }{\Delta_N(\xx )} \det_{i,j=1} ^{N}[p_{\mathfrak{s},\alpha ,t}^{(N)}(x_i,y_j) ]
,\end{align*}
which gives the transition density of non-colliding systems of $N$-particles $L_{\mathfrak{s},\alpha }^{(N)}$-diffusions.

By the same computation for the h-transformation as in \cite{KoO01}, the stochastic differential equation for the $N$-particle $L_{\mathfrak{s},\alpha}^{(N)}$-diffusions h-transformed by $\Delta_N(\mathbf{x})$ is given by \eqref{:41a}.
Hence, the density $\pp_{\mathfrak{s},\alpha ,t}^N$ is the transition density of the solution to the equation \eqref{:41a} once uniqueness of the solution to \eqref{:41a} is established. The required uniqueness is proved in the following lemma.
Thus, for any $\xx \in \intnWN $, we obtain
\begin{align*}
  T_{\mathfrak{s},\alpha,t}^N (\xx,d\yy )=e^{-\lambda_{\mathfrak{s}}^N t} \frac{\Delta_N(\yy ) }{\Delta_N(\xx )} \det_{i,j=1} ^{N}[p_{\mathfrak{s},\alpha ,t}^{(N)}(x_i,y_j) ] d\yy
.\end{align*}

\begin{lemma}\label{l:32}
Let $\mathfrak{s}\in \R$ and $\alpha >-1$.
Then, for any starting point $\xx\in\nWN $, the stochastic differential equation \eqref{:41b} has a unique strong solution, and the solution satisfies the non-explosion and non-colliding property.
\begin{proof}
  We first recall a sufficient condition for strong uniqueness from \cite[Theorem 2.2]{GrM14}.
  Set $I =[0, \infty)$. 
  For continuous functions $\map{\sigma ,b}{I}{\R}$ and a continuous function $\map{H}{I ^2}{[0,\infty)}$, we consider the following stochastic differential equation of $N$-particles on $I$: for $i=1,\ldots, N$, 
  \begin{align*}
  dX_t^i=\sigma(X_t^i)dB_t^i +\Big( b(X_t^i) +\sum_{j\neq i}^{N}\frac{H(X_t^i, X_t^j)}{X_t^i -X_t^j} \Big)dt
  .\end{align*}
  This  has a unique strong solution, which is non-explosive and non-colliding, provided that the following conditions (C1)--(C2), (A1)--(A4) hold on $I$:
  \begin{enumerate}
    \item[(C1)]
    There exists a function $\map{\rho }{(0,\infty )}{(0,\infty)}$ such that $\int _{0^+} \rho ^{-1} (x) dx =\infty $ and that $|\sigma (x) -\sigma (y) |^2 \le \rho (|x-y|)$.
    \hangindent=23pt
     Moreover, $b$ is Lipschitz continuous.

\item[(C2)]
  There exists a constant $c>0$ such that $\sigma (x)^2 +b(x) x \le c (1+ x^2 )$ and $ H(x,y) \le c(1+ xy )$.

  \item[(A1)]
  For $0\le w<x<y<z$, it holds that $ H(w,z) (y-x) \le H(x,y) (z-w)$.

  \item[(A2)]
  There exists a constant $c \ge 0$ such that $\sigma ^2(x) +\sigma ^2(y) \le c(x-y)^2 +4H(x,y).$

  \item[(A3)]  
  There exists a constant $c \ge 0$ such that for any $0\le x< y<z$,
  \begin{align*}
  H(x,y)(y-x) +H(y,z) (z-y) \le c (z-y)(z-x)(y-x) +H(x,z)(z-x)
  .\end{align*}

\item[(A4)]
The set $G :=\{x \st \sigma^2(x)+H(x,x) =0\}$ consists of isolated points.
Furthermore, for any $x\in G$ and  $y_1,\ldots, y_{N-2} \in I $, we have  
\begin{align*}
 b(x) +\sum_{j=1}^{N-2} \frac{H(x, y_j)}{x-y_j } \mathbf{1}_{I \setminus \{x\} } (y_j ) \neq 0
.\end{align*}
\end{enumerate}

We specialise this result to equation \eqref{:41b} by taking 
\begin{align*}
\sigma (x)= \sqrt{2x(1+x)},\quad b(x)=-\mathfrak{s}x+N+\alpha  ,\quad H(x,y)= 2xy +x+y
.\end{align*}
Then, it remains to show (C1)--(C2), (A1)--(A4) to prove this lemma.

Conditions (C1), (C2), and (A2) are immediate. 
Let $f(x,y )=H(x,y)/(y-x)$ for $ x <y$.
A direct calculation shows that $\partial_{x} f(x,y ) \ge 0$ and $\partial_{y} f(x,y ) \le 0$ for all $x,y\ge 0$. 
It follows that $ f(w,z) \le f(x,y)$ for $ w< x <y <z$, which implies that (A1) holds.
The condition (A3) is satisfied for $c\ge 2$ because 
 \begin{align*}
H(x,z)(z-x) -  H(x,y)(y-x) -H(y,z) (z-y) 
 =-2 (z-x) (y-x)(z-y)
.\end{align*}
Since $G=\{0\}$, we see (A4) from $\alpha >-1$.  
Thus, we have shown all conditions, and therefore the proof of this lemma is completed.
\end{proof}
\end{lemma}

\begin{lemma}\label{l:33}
Assume that $\mathfrak{s}\in \R$ and $\alpha >-1$.  
Then, the semigroup $\{ T_{\mathfrak{s},\alpha, t}^{N} \}_{t\ge 0} $ is Feller-Dynkin.
 \begin{proof}
Fix $f \in  C_{\infty} (\nWN )$.
  By a standard argument, the proof of $T_{\mathfrak{s},\alpha,t}^{N} f (\xx ) \in C_{\infty} (\nWN )$ is reduced to showing that, for sufficiently large $R$, 
\begin{align}\label{:A28}
  \limi{\xx } P_{\xx} (\XX_t \in [0,R]^N)= 0
 , \\ \label{:A29}
   \lim_{\xx \to \yy } \int_{\nWN \cap  [0,R]^N } |\mathbf{p}_{\mathfrak{s},\alpha, t}^{N}(\xx,\zz )-\mathbf{p}_{\mathfrak{s},\alpha, t}^{N} (\yy,\zz )  |d\zz =0
.\end{align}

We prove \eqref{:A28} and \eqref{:A29} for the most singular case $\xx=(x,\ldots,x)$ and $\yy =(y,\ldots,y)$.
In this case, we have 
\begin{align*}
\pp_{\mathfrak{s},\alpha ,t}^N(\mathbf{x},\mathbf{y})
=\frac{e^{-\lambda_{\mathfrak{s}}^N t} \Delta_N(\mathbf{y} )}{\prod_{k=1}^{N-1}k! }\det_{i,j=1}^N \big[  \partial_x^{i-1} p_{\mathfrak{s},\alpha, t}^{(N)} (x,y_j) \big]
.\end{align*}
Therefore, for a constant $\Ct{c;43}$ independent of $x,y$, we obtain 
\begin{align*}
  P_{\xx} (\XX_t \in [0,R]^N)
  &\le 
  \Ctref{c;43} \prod_{i=1}^{N} \int_0^R \big|  \partial_x^{i-1} p_{\mathfrak{s},\alpha, t}^{(N)} (x,y) \big| dy
  \\&
  \le  \Ctref{c;43}  \Big(\Ctref{c;A10} \int_0^R J(y) dy\Big)^{N-1} \int_0^R  p_{\mathfrak{s},\alpha, t}^{(N)} (x,y)  dy 
,\end{align*}
from \eqref{:A23}, where $J \in L^1([0,R])$ is given in \ref{s:A1}.
Hence, we have \eqref{:A28} by \eqref{:A22}.
Furthermore, by \lref{l:A2} and the dominated convergence theorem, we obtain \eqref{:A29}.

The general case is treated by decomposing the limiting point into clusters of equal coordinates and rewriting the determinant by divided differences within each cluster; see the proof of \cite[Lemma 4.3]{Ass20}.
By the repeated mean value theorem, these divided differences converge to the corresponding derivatives at the common limiting coordinate. 
Hence, the same computations as above yield \eqref{:A28} and \eqref{:A29} for general $\xx,\yy\in \nWN$.

Since we have proved $T_{\mathfrak{s},\alpha, t}^{N}(C_\infty (\nWN )) \subset C_\infty(\nWN )$ and the process has continuous paths, the strong continuity follows from \cite[III Lemma 6.7]{RoW00}.
\end{proof}
\end{lemma}

\begin{lemma}\label{l:34}
Assume that  $\mathfrak{s} >-1$ and $\alpha>-1$.
Then, the Pickrell ensemble $\msaN $ is the unique invariant probability measure of $\{ T_{\mathfrak{s},\alpha ,t}^{N} \}_{t\ge 0}$.
\begin{proof}

This follows by the same proof as of \cite[Proposition 4.4]{Ass20}.
More precisely, from the fact that $p_{\mathfrak{s},\alpha ,t}^{(N)}$ is symmetric with respect to the speed measure $m_{\mathfrak{s}, \alpha }^{(N)}$ and the equation $ \int_{\nWN} \Delta_N(\xx ) \det_{i,j=1} ^{N}[p_{\mathfrak{s},\alpha ,t}^{(N)}(y_i,x_j) ]  d\xx = e^{\lambda_{\mathfrak{s}}^N t}\Delta_N(\yy )$, we have $ m_{\mathfrak{s}, \alpha }^N T_{\mathfrak{s}, \alpha , t}^N = m_{\mathfrak{s}, \alpha }^N $ by a straightforward computation.
Furthermore, we have $T_{\mathfrak{s}, \alpha ,t}^N(\xx ,A)>0 $ for any $\xx \in \mathring{W}_{\ge}^{N}$ and a Borel set $A \subset \mathring{W}_{\ge}^{N}$ with positive Lebesgue measure from \cite[Theorem 4]{KaM59} . 
Therefore, the uniqueness follows by the same argument as in \cite[Proposition 4.4]{Ass20}.
\end{proof}
\end{lemma}

\subsection{Dual operators and h-transformations}
To establish the intertwining relation of the Pickrell diffusion, we employ the same strategy as in \cite{BuK25}.
Key equations in this technique are two h-transformations, which we prepare in this subsection.

We consider the operator $L_{\mathfrak{s},\alpha}^{(N)}$ for all $\alpha \in\R$.
If $\alpha \le -1$, the origin is an exit boundary.
Let 
\begin{align*}
\hat{L}_{\mathfrak{s},\alpha }^{(N)} 
:=x(1+x) \frac{d ^2}{d x^2}+ \{ (2N+\mathfrak{s} )x -\alpha\} \OD{}{x}
\end{align*}
be the Siegmund dual operator of $L_{\mathfrak{s},\alpha}^{(N)}$ on $[0,\infty)$. 
Here, the point $\infty$ is a natural boundary; the origin is exit for $\alpha \ge 0$, regular absorbing for $-1<\alpha <0$, and entrance for $\alpha \le -1$.
The speed measure for $\hat{L}_{\mathfrak{s},\alpha }^{(N)} $ has density 
\begin{align*}
\hat{m}_{\mathfrak{s},\alpha}^{(N)} (x) 
= x^{-(\alpha+1)}(1+x)^{2N+\mathfrak{s}+\alpha-1}
\end{align*}
on $[0,\infty)$.
Let $ \hat{p}_{\mathfrak{s},\alpha ,t}^{(N)}$ be the transition density associated with $\hat{L}_{\mathfrak{s},\alpha}^{(N)}$ with respect to the Lebesgue measure.
We need two formulas of Doob's h-transform.
\begin{lemma}\label{l:44}
  For any $\mathfrak{s} \in \R $,  the following identities hold for $x,y>0$, where
\begin{align*}
c_{\mathfrak{s}}^{N+1}=-2N-\mathfrak{s}, \qquad d_{\mathfrak{s}, \alpha}^{N+1} =-\alpha(2N+\mathfrak{s}+\alpha-1)
.\end{align*}
  \begin{itemize}
  \item[(i)]
For $\alpha \in\R$, we have
  \begin{align}\label{:22a}
  e^{-c_{\mathfrak{s}}^{N+1} t}\hat {p}_{\mathfrak{s},\alpha ,t}^{(N+1)} (x,y)\frac{(\hat{m}_{\mathfrak{s},\alpha }^{(N+1)} (y))^{-1}}{(\hat{m}_{\mathfrak{s},\alpha}^{(N+1)} (x))^{-1}  }= p_{\mathfrak{s},\alpha +1,t}^{(N)} (x,y)
  .\end{align}
  
  \item[(ii)]
  For  $\alpha >-1$, we have
  \begin{align}\label{:22b}
  e^{-d_{\mathfrak{s},\alpha}^{N+1}  t}p_{\mathfrak{s}+2\alpha -2, -\alpha ,t}^{(N+1)} (x,y) \frac{y^{\alpha }}{x^{\alpha }}= p_{\mathfrak{s},\alpha  ,t}^{(N)} (x,y)
.  \end{align}
    
  \end{itemize}
  \begin{proof}
The proof follows from straightforward computations.
Indeed, we see that $(\hat{m}_{\mathfrak{s},\alpha}^{(N+1 )} (x ))^{-1} =x^{\alpha+1} (1+x)^{-2 N-\mathfrak{s}-\alpha-1}$ is a positive eigenfunction of  $\hat{L}_{\mathfrak{s},\alpha}^{(N+1)} $ with eigenvalue $c_{\mathfrak{s}}^{N+1}$.
 Hence \eqref{:22a} follows from
  \begin{align*}
  \hat{L}_{\mathfrak{s},\alpha }^{(N+1)} \circ (\hat{m}_{\mathfrak{s},\alpha}^{(N+1)})^{-1} 
  =c_\mathfrak{s}^{N+1} (\hat{m}_{\mathfrak{s},\alpha }^{(N+1)})^{-1}  + (\hat{m}_{\mathfrak{s},\alpha }^{(N+1)})^{-1}  L_{\mathfrak{s},\alpha +1}^{(N)} 
  .\end{align*}
 Similarly, since $L_{\mathfrak{s}+2\alpha -2, -\alpha}^{(N+1)} x^{\alpha}= d_{\mathfrak{s},\alpha}^{N+1} x^{\alpha }$, identity \eqref{:22b} follows from 
\begin{align*}
L_{\mathfrak{s}+2\alpha -2 , -\alpha}^{(N+1)} \circ x^{\alpha }
=d_{\mathfrak{s},\alpha }^{N+1} x^{\alpha } +x ^{\alpha} L_{\mathfrak{s}, \alpha}^{(N)} 
.\end{align*}
  \end{proof}
\end{lemma}

Hereafter, to simplify notation, we write the characteristic functions
\begin{align*}
  \mathbf{1}_{N,\xx}^{+} (\yy)=\mathbf{1}_{\WNNp (\xx)} (\yy), \qquad \mathbf{1}_{N,\xx}(\yy)=\mathbf{1}_{\nWNN (\xx)} (\yy)
.\end{align*}
 
\begin{lemma}
For any $\mathfrak{s}\in\R $ and $\alpha>-1$, the following equations hold:
\begin{align}& \label{:56a}
  \int d\zz \, \mathbf{1}_{ N,\zz }^{+}(\yy) \det_{i,j=1}^{N+1}  [ p _{\mathfrak{s},\alpha,t}^{(N+1)}(x_i,z_j)]   
  \\ &\notag  \qquad 
  = \int d\zz \,  \mathbf{1}_{ N,\xx }^{+} (\zz )  \det_{i,j=1}^N  [\hat{p}_{\mathfrak{s},\alpha,t}^{(N+1)}(z_i,y_j)] \prod_{k=1}^N \frac{\hat{m}_{\mathfrak{s},\alpha} ^{(N+1)} (z_k) }{ \hat{m}_{\mathfrak{s},\alpha} ^{(N+1)}  (y_k) } 
,\\ &\label{:56b}
  \int d\zz \,  \mathbf{1}_{ N,\zz }(\yy)  \det_{i,j=1}^{N}  [ p _{\mathfrak{s},\alpha,t}^{(N+1)}(x_i,z_j)] 
  \\ & \notag \qquad
  = \int d\zz \, \mathbf{1}_{ N,\xx } (\zz )  \det_{i,j=1}^N  [\hat{p}_{\mathfrak{s},\alpha,t}^{(N+1)}(z_i,y_j)] \prod_{k=1}^N \frac{ \hat{m}_{\mathfrak{s},\alpha} ^{(N+1)}  (z_k) }{ \hat{m}_{\mathfrak{s},\alpha} ^{(N+1)}  (y_k)}
.\end{align}
\begin{proof}
These are specific formulas of \cite[(13.26)]{AOW19} in our setting.
See also \cite[Remark 4]{BuK25} for a direct computational proof.
\end{proof}
\end{lemma}

\subsection{Intertwining relations and boundary Feller-Dynkin processes}
The proof of \tref{t:35} is reduced establishing to two shifted intertwining relations.
The first of these is the shifted intertwining relation with respect to $\LNNp $, which is proved using the  technique in \cite[Theorem 5.1]{Ass20}.

\begin{lemma}\label{l:23}
Assume that $\mathfrak{s} \in\R $ and $\alpha>-1$.
Then, for any $N\in\N $, $f\in C_{\infty} (W_{\ge}^N)$, and $t \ge 0$, we have
\begin{align}\label{:23a}
T_{\mathfrak{s},\alpha ,t}^{N+1} \LNNp f=\LNNp T_{\mathfrak{s},\alpha +1,t}^{N}f
.\end{align}
\begin{proof}
It is sufficient to show \eqref{:23a} for $t>0$.
We first consider the case where $\xx\in \intnWN $.
Multiplying both sides in \eqref{:56a} by
\begin{align*}
  e^{-\lambda _\mathfrak{s}^{N+1} t} \frac{N! \Delta_N(\yy )}{\Delta _{N+1} (\xx )} d\yy 
,\end{align*}
we see that the left-hand side becomes 
\begin{align*}&
\Big( \int d\zz \,  \pp_{\mathfrak{s},\alpha ,t}^{N+1}(\xx, \zz)  \mathbf{1}_{N, \zz }^{+}(\yy)  \frac{N! \Delta_N(\yy )}{\Delta _{N+1} (\zz )}\Big) d\yy 
=(T_{\mathfrak{s},\alpha ,t}^{N+1} \LNNp  )(\xx,d\yy )
.\end{align*}
On the other hand, using the fact that $\lambda _\mathfrak{s}^{N+1}=N c_\mathfrak{s}^{N+1} + \lambda _\mathfrak{s}^{N}$ and  \eqref{:22a}, we see that the right-hand side becomes
\begin{align*}&
\Big( \int d\zz \, \mathbf{1}_{ N, \xx }^{+} (\zz )  \det_{i,j=1}^N  [\hat{p}_{\mathfrak{s},\alpha ,t}^{(N+1)}(z_i,y_j)]  \prod_{k=1}^N \frac{\hat{m}_{\mathfrak{s},\alpha}^{(N+1)}  (z_k) }{\hat{m}_{\mathfrak{s},\alpha } ^{(N+1)} (y_k)}   \ts e^{-\lambda _\mathfrak{s}^{N+1}t} \frac{N! \Delta_N(\yy )}{\Delta _{N+1} (\xx )} \Big)d\yy 
\\&
=\Big( \int d\zz \, \mathbf{1}_{ N, \xx }^{+}(\zz) \det_{i,j=1}^N [ p _{\mathfrak{s}, \alpha +1, t}^{(N)}(z_i,y_j)] e^{-\lambda_\mathfrak{s}^N t} \frac{N! \Delta _N(\yy )}{\Delta _{N+1} (\xx )} \Big) d\yy  
\\&
=(\LNNp T_{\mathfrak{s}, \alpha +1, t}^{N}  )(\xx,d\yy)
.\end{align*}
Combining these, we obtain 
\begin{align*}
T_{\mathfrak{s},\alpha ,t}^{N+1} \LNNp (\xx, d\yy )=\LNNp T_{\mathfrak{s},\alpha +1,t}^{N}(\xx, d\yy)
,\end{align*}
which implies $T_{\mathfrak{s},\alpha ,t}^{N+1} \LNNp f(\xx )=\LNNp T_{\mathfrak{s},\alpha +1,t}^{N}f (\xx )$ for $\xx \in \intnWN $. 
We can extend this for $\xx \in \nWN $ because of the Feller property \eqref{:31c} and the Feller-Dynkin property established in \lref{l:33}.
Thus we complete the proof.
\end{proof}
\end{lemma}

The second shifted intertwining relation concerns $\LaNN$.
\begin{lemma}\label{l:48}
Assume that $\mathfrak{s} \in \R $ and $\alpha>-1$.
Then, for any $N\in\N$, $f\in C_{\infty} (W_{\ge}^N)$, and $t \ge 0$, we have
  \begin{align*}
  T_{\mathfrak{s}, \alpha +1,t}^{N} \LaNN  f=\LaNN T_{\mathfrak{s},\alpha,t}^{N }f
  .\end{align*}
  \begin{proof}
  Suppose $t>0$ and $\xx\in \intnWN$.
  Then, by \eqref{:32c}, the equality 
\begin{align}\label{:48a}
  T_{\mathfrak{s}, \alpha +1 ,t}^{N } \LaNN (\xx,d\yy)=\LaNN T_{\mathfrak{s},\alpha ,t}^{N}(\xx,d\yy)
\end{align}
 is equivalent to
\begin{align}&\notag
  \int_{ } d\zz \, \mathbf{1}_{ N, \zz }(\yy) \det_{i,j=1}^N[p_{\mathfrak{s}, \alpha+1 ,t}^{(N)}(x_i, z_j)]  \prod_{k=1}^{N} \frac{x_k^{\alpha +1} }{z_k^{\alpha+1}} 
  \\&\label{:32b} \qquad
  =\int_{} d\zz \, \mathbf{1}_{N, \xx }(\zz) \det_{i,j=1}^N  [p_{\mathfrak{s},\alpha ,t}^{(N)}(z_i, y_j)]\prod_{k=1}^{N} \frac{z_k^{\alpha } }{y_k^{\alpha}} 
.\end{align}
Note that $L_{\mathfrak{s}+2\alpha ,-\alpha}^{(N)} = L_{\mathfrak{s}+2(\alpha-1),-\alpha}^{(N+1)}$ by definition, and hence we have
\begin{align}\label{:48b}
  p_{\mathfrak{s}+2\alpha,-\alpha ,t}^{(N) }(z,y)=p_{\mathfrak{s}+2(\alpha-1),-\alpha ,t}^{(N+1) }(z,y)
.\end{align} 
Furthermore, we remark that 
\begin{align}\label{:48c}
  d_{\mathfrak{s},\alpha +1}^{N+1}-c_{\mathfrak{s}+2\alpha}^{N+1} =d_{\mathfrak{s},\alpha }^{N+1} 
.\end{align}
  We obtain \eqref{:32b} by calculation
  \begin{align*}&
  \int_{ } d\zz \, \mathbf{1}_{N,\zz }(\yy) \det_{i,j=1}^N[p_{\mathfrak{s}, \alpha +1 ,t}^{(N)}(x_i, z_j)]  \prod_{k=1}^{N} \frac{x_k^{\alpha +1} }{z_k^{\alpha +1}}  
  \\&
  = e^{-N d_{\mathfrak{s},\alpha +1 }^{N+1} t}  \int_{ } d\zz \, \mathbf{1}_{N, \zz }(\yy)\det_{i,j=1}^{N} [p_{\mathfrak{s}+2\alpha , -\alpha-1 ,t}^{(N+1)}(x_i, z_j)]   && \qquad\text{ from \eqref{:22b}}
  \\&
  = e^{-N d_{\mathfrak{s},\alpha +1}^{N+1} t} \int_{ } d\zz \, \mathbf{1}_{N, \xx }(\zz) \det_{i,j=1}^N  [\hat{p}_{\mathfrak{s}+2\alpha,-\alpha-1 ,t}^{(N+1)}(z_i,y_j)] 
  \\& \qquad \qquad \qquad \qquad 
  \ts \prod_{k=1}^N \frac{ \hat{m}_{\mathfrak{s}+2\alpha, -\alpha-1 }^{(N+1)}   (z_k)}{\hat{m}_{\mathfrak{s}+2\alpha ,-\alpha-1 }^{(N+1)}   (y_k)}  &&  \qquad\text{ from \eqref{:56b}}
  \\&
  =e^{-N (d_{\mathfrak{s},\alpha +1}^{N+1}-c_{\mathfrak{s}+2\alpha}^{N+1}) t}\int_{ } d\zz \, \mathbf{1}_{N, \xx }(\zz)  \det_{i,j=1}^{N} [p_{\mathfrak{s}+2\alpha,-\alpha ,t}^{(N) }(z_i,y_j)] 
  && \qquad\text{ from  \eqref{:22a}  }
 \\& 
 =e^{-N d_{\mathfrak{s},\alpha }^{N+1} t} \int_{ } d\zz \, \mathbf{1}_{N, \xx }(\zz)  \det_{i,j=1}^{N} [p_{\mathfrak{s}+2(\alpha-1),-\alpha ,t}^{(N+1) }(z_i,y_j)] 
 && \qquad\text{ from \eqref{:48b}, \eqref{:48c}}
\\&
 = \int_{ } d\zz \, \mathbf{1}_{N, \xx }(\zz) \det_{i,j=1}^N [p_{\mathfrak{s},\alpha ,t}^{(N)}(z_i,y_j)]  \prod_{k=1}^{N} \frac{z_k^{\alpha}  }{y_k^\alpha}&& \qquad\text{ from \eqref{:22b}}
.\end{align*}
Thus, we have proved \eqref{:48a} for $\xx\in \intnWN$.
From the Feller property \eqref{:31d} and the Feller-Dynkin property established in \lref{l:33}, we can extend \eqref{:48a} to all $\xx\in \nWN$, which completes the proof.
  \end{proof}
  \end{lemma}

 \medskip
 \noindent
{\em Proof of \tref{t:35}  } \quad 
   From \lref{l:23} and \lref{l:48}, we have the identity $ T_{\mathfrak{s},\alpha, t}^{N+1} \LNNp \LaNN = \LNNp \LaNN T_{\mathfrak{s},\alpha, t}^{N}  $. 
  This equation with \lref{l:35} concludes the statement of this theorem.
\qed

\begin{corollary}
For any $\mathfrak{s}>-1$ and $\alpha >-1$, we have
\begin{align} \label{:49b}
m_{\mathfrak{s},\alpha }^{N+1} \LNNp =m_{\mathfrak{s},\alpha+1 }^{N}
,\\ \label{:49c}
m_{\mathfrak{s},\alpha +1}^{N} \LaNN =m_{\mathfrak{s},\alpha }^{N}
.\end{align}
Furthermore, combining these, we have
\begin{align}\label{:49a}
  m_{\mathfrak{s},\alpha }^{N+1} \LaNNp =m_{\mathfrak{s},\alpha }^{N}
.\end{align}
\begin{proof}
Equations \eqref{:49b} and \eqref{:49c} follow from \lref{l:34} with \lref{l:23} and \lref{l:48}, respectively.
\end{proof}
\end{corollary}

\begin{remark}\label{r:410}
  \begin{itemize}
  \item[(i)]
For $\alpha \in \{0\} \cup \N$, the equation \eqref{:49a} immediately follows from \eqref{:31h} since $m_{\mathfrak{s},\alpha }^{N}$ is the distribution of the radial parts of a $\mathbb{U}(N+\alpha) \ts \mathbb{U}(N)$-invariant random matrix determined by \eqref{:13a}.

  \item[(ii)]
For $a ,b >-1$, let $m_{\mathrm{Jac}, a , b }^N$ be the Jacobi ensemble on $[0,1]$ given by
\begin{align*}
  m_{\mathrm{Jac}, a , b }^N (d\uu )=\frac{1}{Z_{\mathrm{Jac}, a ,b }^N} \Delta_{N} ^2(\uu ) \prod_{k=1}^{N} u_k^{a } (1-u_k)^{b } d\uu 
,\end{align*}  
where $Z_{\mathrm{Jac}, a ,b }^N$ is the normalising constant.
The change of variables 
\begin{align*}
  x_i=\frac{u_i}{1-u_i }, \qquad y_i=\frac{v_i}{1-v_i }
\end{align*}
transforms the equations \eqref{:49b} and \eqref{:49c} into the identities 
\begin{align}\label{:410a}
m_{\mathrm{Jac}, \alpha, \mathfrak{s}}^{N+1} L_{ N}^{N+1}  = m_{\mathrm{Jac}, \alpha+1, \mathfrak{s}}^{N}
,\\ \label{:410b}
m_{\mathrm{Jac}, \alpha+1, \mathfrak{s}}^{N} L_{\alpha, N}^N=m_{\mathrm{Jac}, \alpha, \mathfrak{s}}^{N}
.\end{align}
Here, the kernels $  L_{ N}^{N+1} : \WNp \cap [0,1]^{N+1} \dashrightarrow \WN \cap [0,1]^{N}$ and $  L_{\alpha, N}^{N} : \WN \cap [0,1]^{N} \dashrightarrow \WN \cap [0,1]^{N}$ are given by
\begin{align*}
  L_{ N}^{N+1} (\uu, d\vv )&:= N! \frac{\Delta_N(\vv)}{\Delta_{N+1} (\uu)}  \frac{\prod_{k=1}^{N+1}(1-u_k)^{N}}{\prod_{k=1}^N (1-v_k)^{N+1}} \mathbf{1}_{W_{}^{N,N+1}(\uu )}(\vv) d\vv 
,  \\
  L_{\alpha, N}^N (\uu, d\vv )&:= (\alpha+1)_N \frac{\Delta_N(\vv)}{\Delta_N(\uu)} \prod_{k=1}^N \frac{v_k^\alpha }{u_k^{\alpha +1}}\frac{(1-u_k)^{N+\alpha}}{(1-v_k)^{N+\alpha+1}} \mathbf{1}_{W_{\ge}^{N,N}(\uu )}(\vv) d\vv 
.\end{align*}
Equation \eqref{:410a} corresponds to the case $\beta = 2$ of the $\beta$-Jacobi corners process introduced in \cite[Definition 2.6]{BoG15}; this process is defined for general $\beta$.
It would be interesting to find a general $\beta$ analogue for \eqref{:410b}.
\end{itemize}
\end{remark}

\medskip
\noindent
{\em Proof of \tref{t:12}  } \quad
We have checked \eqref{:39a} and \eqref{:39b} from \tref{t:35} and \lref{l:48}, respectively.
Therefore, the first assertion immediately follows from \pref{p:39}.
Furthermore, the existence of unique invariant probability measure results from \lref{l:24} with \lref{l:34}.
\qed

\section{Approximation of boundary Feller-Dynkin processes}\label{s:5}
\subsection{A general approximation result}
Define an embedding map $\map{\mathfrak{r}_N }{\nWN }{\Omega}$ by $\mathfrak{r}_N (\xx )=(\boldsymbol{a }(\xx), \gamma(\xx ))$, where
\begin{align*}
a_i(\xx ) =
\begin{cases}
\frac{x_{N+1-i} }{N^2} &\text{ for } i\le N,
\\
0 & \text{ for }  i\ge N+1
,\end{cases}
\qquad
\gamma (\xx )=\frac{x_1+\ldots+x_N}{N^2}
\end{align*} 
for $\xx =(x_1\ldots, x_N) \in \nWN $.

The boundary probability measure $m \in \MMp (\Omega ) $ is approximated by the corresponding coherent measures in the following sense:
\begin{lemma}\label{l:51}\cite{BoO01, Buf15}
For $\{m^N\}_{N\in\N } \in \varprojlim \mathcal{M}_{p,\alpha } (\nWN ) $, let $m \in \MMp (\Omega ) $ be the corresponding measure.
Then, we have
\begin{align*}
  \limi{N} (\mathfrak{r}_N)_* m^N =m  \quad \text{ in distribution}
.\end{align*}
\end{lemma}
As a dynamical version of this result, we obtain the following:
\begin{lemma}\label{l:52}
For $\{m^N\}_{N\in\N } \in \varprojlim \mathcal{M}_{p,\alpha } (\nWN ) $, let $m \in \MMp (\Omega ) $ be the corresponding measure.
Assume that a family of Feller-Dynkin processes $\{ \XX ^N \}_{N\in\N }$ is coherent with respect to $\{ \nWN, \LaNNp \} _{N\in\N}$, and let $\XX ^\infty$ be its boundary Feller-Dynkin process on $\Omega $. 
Assume that $\XX _0^N$ has distribution $m^N$,  and $\XX _0 ^\infty $ has distribution $ m $. 
Then, for any $t\ge 0$ fixed, we have
\begin{align*}
\limi{N}\mathfrak{r}_N(\XX _t ^N) =\XX _t^\infty \quad \text{ in distribution}
.\end{align*}
\begin{proof}
This statement can be proved by the same argument as in \cite[Proposition 5.4]{Ass20}.
From the equality
\begin{align*}
  m^{N+1} T_{t}^{N+1} \LaNNp = m^{N+1}  \LaNNp T_{t}^N =m^N T_{t}^N
,\end{align*}
it follows that $\{ m^N T_{t}^N \}_{N\in\N} \in \varprojlim \mathcal{M}_{p,\alpha } (\nWN ) $. 
Furthermore, the corresponding measure of this coherent family is $m T_{t}^\Omega$.
Therefore, by applying \lref{l:51}, we conclude the proof.
\end{proof}
\end{lemma}

We note that stronger approximations, such as those discussed in \cite[Section 2]{AsM}, may hold, but the above lemma is sufficient for our purpose in this paper.

\subsection{Boundary Laguerre processes}
We first deal with the Laguerre process (without the Ornstein--Uhlenbeck term).
We consider the $N$-dimensional stochastic differential equation
\begin{align} \label{:51b}
  dY_t^{N, i} &=\sqrt{2 Y_t^{N,i} } dB_t^i +\Big( \alpha +N +\sum_{j\neq i}^{N} \frac{ Y_t^{N,i}+Y_t^{N,j} }{Y_t^{N,i} -Y_t^{N,j} } \Big)dt
\end{align}
for $i=1,\ldots, N$.
Similar to \eqref{:15a}, equation \eqref{:51b} has a unique strong solution in $\nWN $ for any starting point $\yy \in \nWN$ for $\alpha >-1$ \cite[Corollary 6.6]{GrM14}.
We call the unique strong solution $\YY ^N =(Y^{N,i})_{i=1}^N$ the Laguerre process of parameter $\alpha$, and let $ \mathsf{T}_{\alpha}^N := \{ \mathsf{T}_{\alpha, t}^N \}_{t\ge 0}$ be the Markov semigroup, which is  Feller-Dynkin \cite[Proposition 1.3]{Ass19b}.
There exists a boundary process associated with the family of the Laguerre processes.

\begin{theorem}\label{t:51}
There exists a unique Feller-Dynkin semigroup $\{ \mathsf{T}_{t}^{\Omega } \}_{t\ge 0}$ on $\Omega $ such that, for any $\alpha \in \{0\} \cup \N$, $N\in\N$, $f\in C_{\infty} (W_{\ge}^N)$, and $t \ge 0$, we have
\begin{align*}
\mathsf{T}_{t}^{\Omega } \LaON f=\LaON \mathsf{T}_{\alpha ,t}^{N } f 
.\end{align*}
Furthermore, the boundary process $(\boldsymbol{a } (t), \gamma (t))$ associated with $\{ T_{t}^{\Omega }  \}_{t\ge 0}$ is given by 
\begin{align}\label{:51c}
     a_i (t)=a_i(0), \quad \gamma(t)=t+\gamma(0)
    .\end{align}
\begin{proof}
Since the family $ \{ \mathsf{T}_{\alpha}^N \}_{N\in\N}$ is coherent with respect to $\{ \nWN , \LaNNp \}_{N\in\N}$ \cite[Theorem 1.7]{KaS26} and \eqref{:39b} holds \cite[Theorem 1.8]{KaS26}, we obtain $\{\mathsf{T}_{t}^{\Omega }\}_{t\ge 0}$ by \pref{p:53}.

Let $\YY_0^N \in\nWN $ be a random variable with distribution $\delta_{(\boldsymbol{a } (0), \gamma (0))} \LaNNp $.
Let $A^N \in M_{N+\alpha, N}(\C)$ be a random matrix whose squared singular values are $\YY_0^N $. 
Let $W_t^N$ be the $M_{(N+\alpha)\ts N}(\C) $-valued Brownian motion starting at zero, independent of $A^N$.
Then, the ordered eigenvalue process of $ (A^N+W_t^N)^*(A^N+W_t^N) $ has the same law as the Laguerre process.
Since only distributional properties of $\YY^N$ are relevant, we may and do
realise $\YY^N $ as the eigenvalue process.
Then, Weyl's inequality yields
\begin{align}
|a_i(\YY_t^N ) -a_i(\YY_0^N) |
  \label{:51d}
\le \frac{1}{N^2} \{ \| (W_t^N)^* W_t^N \| + 2 \| (A^N)^*W_t^N \| \}
,\end{align}
where $\| \cdot \|$ denote the operator norm, which is equal to the largest singular value.

From $\limi{N} \sup_{0\le t\le T} \| (W_t^N)^* W_t^N \| /N=4 T$ a.s.~proved in \cite[Theorem 4.1]{PeT09}, we obtain $\sup_{0\le t\le T}\|  W_t^N \| =O ( \sqrt{N})$ a.s.
Furthermore, we have $ \| (A^N)^*W_t^N \| \le \sqrt{Y_0^{N,N}}  \| W_t^N \| $, and $\sqrt{Y_0^{N,N}}/N $ is bounded in probability from \lref{l:51}.
Hence, \eqref{:51d} implies that $\max_{1\le i\le N}\sup_{0\le t\le T}|a_i(\YY_t^N ) -a_i(\YY_0^N) |$ converges to zero as $N\to\infty $ in probability.
Therefore, by  \lref{l:52},  we obtain $a_i (t)=a_i(0)$.

To obtain the equation for $\gamma$, by a straightforward computation, we have
  \begin{align*} 
  d \gamma (\YY_t^N) &=\frac{\sqrt{2\gamma (\YY _t^N)} }{N} d B_t +\frac{N+\alpha }{N}dt
,\end{align*}
where $ B_t$ is the standard Brownian motion.
Letting $N\to \infty$ we have $\partial_t \gamma(t)=1$, which yields $\gamma(t)=t+\gamma(0)$.
\end{proof}
\end{theorem}

\medskip
\noindent
{\em Proof of \tref{t:54}  } \quad
The existence of the boundary Feller-Dynkin process $\{T_{\alpha, t}^{\Omega }\}_{t\ge 0}$  results from \pref{p:39} with \pref{p:53}.
Here, the $\alpha$-independence follows from the shifted intertwining relation $  T_{\alpha+1,t}^{N} \LaNN  =\LaNN T_{\alpha,t}^{N}$, that was proved in \cite[Lemma 14]{BuK25}.

By a deterministic time change and It\^o's formula, we have $\{ X_t^{N,i}\}_{t\ge 0}\overset{\mathrm{law}}{=}\{ e^{-t} Y_{e^ t-1}^{N,i} \}_{t\ge 0}$, which yields $a_i (\XX_t^{N} ) =e^{-t} a_i( \YY_{e^ t-1}^{N} )$ and $\gamma (\XX_t^{N} ) =e^{-t} \gamma ( \YY_{e^ t-1}^{N} )$.
Therefore, equation \eqref{:54a} follows from \tref{t:51}.
\qed

\section{Discrete kernels that converge to  $\LaNNp $ }\label{s:6}
In this section, we show that our kernel $\LaNNp $ is obtained from the branching formula for  multivariate Jacobi polynomials.
\subsection{Branching formula for multivariate Jacobi polynomials}
Fix two parameters $\alpha, \beta > -1$.
Let $\mathfrak{p} _{n} (x ;\alpha ,\beta )  $  be the Jacobi polynomials, which are expressed in terms of the Gauss hypergeometric function ${}_2 F_{1}$ as 
\begin{align*}
  \mathfrak{p} _{n} (x ;\alpha ,\beta ) =\frac{\Gamma(n+\alpha +1 )}{\Gamma(n+1 )\Gamma(\alpha +1 )} {}_2 F_{1}
\Big(
  \begin{matrix}
  -n, n+\alpha +\beta +1
  \\
  \alpha +1
\end{matrix}
; \frac{1-x}{2} \Big)
.\end{align*}
These polynomials are orthogonal on $[-1, 1]$ with the weight function $(1-x) ^\alpha  (1+x)^\beta$.
We set $\sigma =(\alpha +\beta+1)/2$.
The value at $x=1$ is
\begin{align*}
  \mathfrak{p} _{n} (1; \alpha ,\beta ) = \frac{\Gamma(n+\alpha +1) }{\Gamma( n+1) \Gamma(\alpha +1)}
\end{align*}
and the leading coefficient in $\mathfrak{p}_{n}(x; \alpha ,\beta ) $ is 
\begin{align*}
  k_{n} := 2^{-n}\frac{\Gamma (2n+2\sigma) }{\Gamma (n+ 2\sigma) \Gamma (n+1)  }
.\end{align*}

Let $\lambda$ be an integer partition and let $l(\lambda ) $ be its length.
Define the multivariate Jacobi polynomial indexed by a partition $\lambda $ with $l(\lambda )\le n$ by 
\begin{align*}
  \mathfrak{P}_{\lambda}(x_1,\ldots, x_n; \alpha , \beta  ) := \frac{\det_{i,j=1}^{n} [\mathfrak{p}_{\lambda_i +n-i} (x_j ; \alpha ,\beta ) ]}{\prod_{1\le i <j\le n }(x_i-x_j)}
.\end{align*}
Let  $1_n $ denote the $n$-dimensional vector whose components are all equal to $1$.

\begin{lemma}\cite[Proposition 7.1]{OlO12}\label{l:61}
  \begin{align*}
    \mathfrak{P}_{\lambda}(1_{n} ; \alpha ,\beta ) 
    &=2^{ - \frac{n(n-1)}{2} } \prod_{1\le i<j\le n} (\lambda_i-\lambda_j +j-i)  
    (\lambda_i+\lambda_j +2n -i-j  +2\sigma )
    \\&\qquad \qquad \quad
    \ts \prod_{i=1}^{n}  \frac{\Gamma(\lambda_i +n-i +\alpha +1)  }{ \Gamma (\lambda_i+n-i+1)  \Gamma(n-i+\alpha +1) \Gamma (i)}
.\end{align*}
\end{lemma}

The Jacobi polynomial $ \mathfrak{P}_{\lambda}$ satisfies a two-step branching formula.
\begin{lemma} \label{l:62}
For any partition $\lambda $ with $l (\lambda) \le n$, we have 
    \begin{align}\label{:62a}
    \mathfrak{P}_{\lambda} (x_1,\ldots, x_{n-1}, 1 ; \alpha ,\beta )= \sum_{\mu \prec \lambda  } \sum_{\nu \prec \mu \cup 0} \frac{c_{\nu,n-1}}{c_{\lambda ,n} }  A_{\mu,\nu }  \mathfrak{P}_\nu (x_1,\ldots, x_{n-1} ;\alpha ,\beta )
  .\end{align}    
Here, 
\begin{align*}
   A_{\mu,\nu }&:=\prod_{i=1}^{n-1}B(\mu_i +n-i-1, \nu _i+n-i-1)
  , \\
   B(m,l)& :=\frac{(2m+\alpha +\beta )\Gamma(m+\beta +1) m! (2l+\alpha +\beta +1) }{2 \Gamma(m+\alpha +\beta +2)\Gamma(m+\alpha +2)}
   \\& \qquad \qquad\qquad \qquad \qquad \qquad\qquad
   \times \frac{\Gamma(l+\alpha +\beta +1) \Gamma(l+\alpha +1)}{\Gamma(l +\beta +1) l!}
  , \\
   c_{\lambda, n} & :=\Gamma(\alpha  + 1)^n \prod_{i=1}^{n} \frac{\Gamma(\lambda_i+n-i+1 )}{\Gamma(\lambda_i+n -i +\alpha +1) }
.\end{align*}
Furthermore, $\mu \prec \lambda $ means 
\begin{align*}
  \lambda_1 \ge \mu_1 \ge \lambda_2 \ge \cdots \ge \lambda_{n-1} \ge \mu_{n-1} \ge \lambda_n
,\end{align*}
and  $\nu \prec \mu \cup 0 $ means 
\begin{align*}
  \mu_1 \ge \nu_1 \ge \mu_2 \ge \cdots \ge \mu_{n-1} \ge \nu_{n-1} \ge 0
.\end{align*}
\begin{proof}
   The proof is directly adapted from that of \cite[Proposition 7.5]{OkO06}.
\end{proof}
\end{lemma}

\subsection{Convergence of discrete kernels}
We set $W_{d,\ge}^{N} := \{\lambda=(\lambda_i)_{i=1}^{N} \in ( \{0\} \cup \N )^N \st \lambda _1 \le \cdots \le \lambda_N\}$.
Define a kernel $L_{N,\alpha, \beta }^{N+1}  :W_{d,\ge}^{N+1} \dashrightarrow W_{d,\ge}^{N}$ by, for  $\lambda \in W_{d,\ge}^{N+1}$, 
\begin{align*}
L_{N,\alpha, \beta }^{N+1} (\lambda , \nu )= \frac{c_{\nu,N}}{c_{\lambda ,N+1} } \frac{\mathfrak{P}_{\nu }(1_{N} ;\alpha ,\beta ) }{\mathfrak{P}_{\lambda}(1_{N+1} ;\alpha ,\beta ) } \sum_{\mu \st \nu \prec \mu \cup 0, \mu \prec \lambda  }   A_{\mu,\nu } 
.\end{align*}
Here, for $\nu =(\nu_1 \le \cdots \le \nu_N) \in W_{d,\ge}^{N}$, we interpret $c_{\nu, N}$ as $c_{(\nu_N,\ldots , \nu_1) , N}$, and other symbols are understood in the same way.
Since $\sum_{\nu\in W_{d,\ge}^{N}} L_{N,\alpha,\beta}^{N+1}(\lambda, \nu )=1$ from \eqref{:62a}, the kernel $L_{N,\alpha, \beta }^{N+1}$ is Markov.
Note that the kernel $L_{N,\alpha, \beta }^{N+1}$ was introduced in
\cite[Section 3.3]{Cue18}.

\begin{lemma}\label{l:63}
For $\lambda \in W_{d,\ge}^{N+1 }$, let $Z_{\lambda } $ be a random variable distributed as $L_{N,\alpha, \beta }^{N+1} (\lambda, \cdot)$.
Then, the distribution of $\kappa ^{-1} Z_{\kappa \lambda }$ converges weakly to the probability measure given by
\begin{align}\label{:63e}
  2^{2N} N! (\alpha+1)_{N}  \frac{\Delta _{N} ( \nu^2 )}{\Delta _{N+1} (\lambda^2 _{})} \bigg( \int_{W_{\ge}^{N+1,N} (\lambda )}^{ }\mathbf{1}_{W_{\ge}^{N,N}(\mu)}(\nu) \prod_{k=1}^{N} \frac{\nu_k^{2\alpha +1 } }{\mu_k ^{2\alpha +1}} d \mu  \bigg)  d\nu 
.\end{align}
Here, for $\lambda =(\lambda_1,\ldots, \lambda_{N+1})$, we write $\lambda^2 =(\lambda_1^2 ,\ldots, \lambda_{N+1}^2)$.
\begin{proof}
 
From $\Gamma (z ) =\sqrt{2\pi }z^{z-1/2} e^z ( 1 +o (1) ) $ as $z\to \infty $, we obtain 
\begin{align} \label{:63b}
  \frac{c_{\kappa \nu,N }}{c_{\kappa \lambda ,N+1} } &=  \frac{\kappa^{\alpha } }{\Gamma(\alpha +1)}\frac{\prod _{i=1}^{N+1 } \lambda _i^{\alpha }}{\prod _{i=1}^{N} \nu_i^{\alpha } } \big(1+ o(1)  \big) 
,\\  \label{:63c}
  A_{\kappa \mu, \kappa \nu} &=2^N \prod_{i=1}^{N} \frac{ \nu_i^{2\alpha +1}} {\mu_i^{2\alpha +1}} \big(1+ o(1)  \big)
\end{align}
for $\kappa \to \infty$.
Because \lref{l:61} yields
\begin{align*}
\mathfrak{P}_{\kappa \lambda}(1_{N+1} ;\alpha ,\beta ) = \frac{\kappa^{N(N+1) +(N+1)\alpha  } \Delta_{N+1}  (\lambda^2 ) \prod_{i=1}^{N+1} \lambda_i^\alpha }{2^{\frac{N(N+1)}{2}} \prod_{i=1}^{N+1} \Gamma(N-i+\alpha +2) \Gamma(i) } \big(1+o(1)  \big)
,\end{align*}
we have
\begin{align} \label{:63d}
  \frac{\mathfrak{P}_{\kappa \nu }(1_{N} ;\alpha ,\beta ) }{\mathfrak{P}_{\kappa \lambda}(1_{N+1} ;\alpha ,\beta )}=\frac{2^{N }  N! \Gamma(N+1+\alpha  ) }{\kappa^{2N +\alpha  } }  \frac{  \Delta_{N} (\nu^2 ) \prod_{i=1}^{N} \nu_i^\alpha }{  \Delta_ {N+1}  (\lambda^2) \prod_{i=1}^{N+1 } \lambda_i^\alpha } \big(1+ o(1)  \big)
.\end{align}

Note that $o(1) $ in \eqref{:63b}, \eqref{:63c}, and \eqref{:63d} are uniform for $\nu $ and $\mu $ in compact subsets.
Hence, combining \eqref{:63b}, \eqref{:63c}, and \eqref{:63d}, we obtain
\begin{align} &\label{:63a}
L_{N,\alpha, \beta }^{N+1}(\kappa \lambda, \kappa \nu ) 
\\\notag
&=\frac{2^{2N} N! (\alpha+1)_{N} }{\kappa ^N} \frac{\Delta _{N} ( \nu^2 )}{\Delta _{N+1} (\lambda^2 _{})} \int_{W_{\ge}^{N+1,N} (\lambda )}^{ } \mathbf{1}_{W_{\ge}^{N,N}(\mu)}(\nu)\prod_{i=1}^{N} \frac{\nu_i^{2\alpha +1 } }{\mu_i ^{2\alpha +1}} d \mu  (1+o(1) )
.\end{align}
Since the error term is uniform on compact subsets, summing \eqref{:63a} against a bounded continuous test function and letting $\kappa\to\infty$, the corresponding Riemann sums converge to the measure in \eqref{:63e}.
This completes the proof.
\end{proof}
\end{lemma}

By the change of variables $\lambda_i^2 =x_i $, $\nu _i^2 =y_i$, and $\mu _i^2= z_i$, the right-hand side of \eqref{:63e} becomes
\begin{align*}
 N! (\alpha+1)_{N}  \frac{\Delta _{N} ( \yy)}{\Delta _{N+1 } (\xx _{})} \bigg( \int_{W_{\ge}^{N+1,N} (\xx )}^{ } \mathbf{1}_{W_{\ge}^{N,N}(\zz )}(\yy ) \prod_{k=1}^{N} \frac{y_k^{\alpha  } }{z_k ^{\alpha +1}} d \zz  \bigg) d\yy  
,\end{align*}
which is identical to $\LaNNp (\xx,d\yy )$ from \lref{l:35} with \eqref{:31g} and \eqref{:32c}.
This naturally leads the following questions.

The boundary processes associated with the kernels $\{ L_{N,\alpha, \beta }^{N+1}\}_{N\in\N } $ were constructed in \cite{Cue18}.
These processes are included in the general framework of non-intersecting birth-and-death chain intertwinings 
developed in \cite{Ass18}.
In view of  \lref{l:63}, it is natural to ask how these boundary processes are related to the boundary Pickrell processes constructed in \tref{t:12}.
In particular, it would be interesting to understand whether this relation can be realised through a suitable diffusive scaling limit.

\appendix
\renewcommand{\thetheorem}{\Alph{section}.\arabic{theorem}}
\section{Estimates for one-particle transition densities }\label{s:A1}
We fix $\alpha>-1$ throughout this section.
Let $p_{t}^{0}(x,y)$ be the transition density associated with the generator
\begin{align*}
L_0 := \frac{1}{2}\frac{d^2}{dx^2}+ \frac{\alpha+\frac{1}{2}}{x}\frac{d}{dx}
\end{align*}
with respect to the speed measure $m_0(dx)=x^{2\alpha+1}dx$.
At the origin, we impose the boundary condition corresponding to that of $X$ introduced in \sref{s:41}: entrance for $\alpha\ge 0$, and regular reflecting for $-1<\alpha<0$.
Let $I_{\alpha }$ denote the modified Bessel function of the first kind of order $\alpha $.
Then,
\begin{align}\label{:A1}
  p_{t}^{0}(x,y)= \frac{1}{t}(xy)^{-\alpha } \exp \Big(-\frac{x^2+y^2}{2t}\Big) I_{\alpha } \Big(\frac{xy}{t}\Big)
.\end{align} 

\begin{lemma}\label{l:A1}
There exist constants $\Ct{c;A1}, \Ct{c;A2}, \Ct{c;A3}, \Ct{c;A4} >0$ such that, for $t>0$ and $x,y\geq0$,
\begin{align}& \notag 
 \frac{ \Ctref{c;A1}}{\sqrt{t} (x+\sqrt{t})^{2\alpha+1}}\exp \Big(- \Ctref{c;A2} \frac{(x-y)^2}{t} \Big)
  \le 
  p_{t}^{0}(x,y)
  \\&\label{:A2}\qquad \qquad\qquad\qquad \qquad
  \le
    \frac{\Ctref{c;A3}}{\sqrt{t} (x+\sqrt{t})^{2\alpha+1}}\exp \Big(- \Ctref{c;A4} \frac{(x-y)^2}{t} \Big)
.\end{align}
Furthermore, there exist constants $\Ct{c;A15}, \Ct{c;A16} >0$ such that, for $k=1,2$,  $t>0$, and $x, y\ge 0$,
\begin{align}\label{:A3}
\Big| \Big ( \frac{\partial }{\partial x } \Big)^k p_{t}^{0}(x,y) \Big|
  &
  \le \frac{ \Ctref{c;A15} }{t^{(k+1)/2} (x+\sqrt{t})^{2\alpha+1}}\exp \Big(-  \Ctref{c;A16} \frac{(x-y)^2}{t} \Big)
.\end{align}
Here, the derivatives at $x=0$ are understood as right derivatives.
\begin{proof}

By the standard asymptotic formulae for $I_\alpha$, namely
\begin{align*}
I_\alpha(z) &\sim \frac{1}{\Gamma(\alpha+1)}\Big(\frac{ z}{2}\Big)^\alpha
 && (z\downarrow0),
\\ 
I_\alpha(z) &= \frac{e^z}{\sqrt{2\pi z}}\Big(1-\frac{4\alpha^2 -1}{8z}+\mathcal{O}(z^{-2}) \Big)
 &&(z\to\infty)
,\end{align*}
there exists a positive constant $\Ct{c;A5} $ such that, for $z>0$,
\begin{align}\label{:A4}
\Ctref{c;A5}^{-1}  z^\alpha(1+z)^{-\alpha-\frac12}e^z
\le 
I_\alpha(z)
\le
\Ctref{c;A5}  z^\alpha(1+z)^{-\alpha-\frac12}e^z
,\\ \label{:A25}
|I_{\alpha+1}(z)-I_\alpha(z) | \le \Ctref{c;A5} z^\alpha (1+z)^{-\alpha -3/2}e^z
.\end{align}

Set $\xi :=x/ \sqrt {t}, \eta :=y/\sqrt{ t}$.
Then, combining \eqref{:A4} with \eqref{:A1}, we obtain
\begin{align}\label{:A12}
\frac{\Ctref{c;A5}^{-1}}{t^{\alpha+1} \left(1+ \xi \eta \right)^{\alpha+\frac{1}{2}}}
  e^{-\frac{(\xi-\eta )^2}{2}}
\le & 
p_{t}^{0}(x,y)
\le
\frac{\Ctref{c;A5}}{t^{\alpha+1} \left(1+\xi\eta \right)^{\alpha+\frac{1}{2}}}  e^{-\frac{(\xi-\eta )^2}{2} }.
\end{align}
A direct computation yields that 
there exist constants $c ,m>0$ such that 
\begin{align}\notag
c^{-1} (1+|\xi-\eta|)^{-m}(1+\xi )^{-2\alpha-1}
\le
(1+\xi\eta )^{-\alpha-\frac{1}{2}}
\le
c (1+|\xi -\eta |)^{m}(1+\xi )^{-2\alpha-1}
.\end{align}
Combining this with \eqref{:A12}, and absorbing the polynomial factors in $1+|\xi-\eta|$ into the Gaussian term, we obtain \eqref{:A2}.

Using $\partial_z (z^{-\alpha } I_\alpha(z)) =z^{-\alpha} I_{\alpha+1}(z)$ and $I_{\alpha +2}=I_\alpha (z)-\{2(\alpha+1)/z\} I_{\alpha+1}(z)$, we have, for $x>0$,
\begin{align*}
&
\partial_x p_{t}^{0}(x,y)
=t^{-\alpha-3/2}(\xi\eta )^{-\alpha} e^{ -\frac{\xi^2+\eta ^2}{2}} \{(\eta-\xi )I_\alpha(\xi\eta )+\eta (I_{\alpha+1}(\xi\eta )-I_\alpha(\xi\eta ) ) \}
,\\&
\partial_x^2 p_{t}^{0}(x,y)
  =t^{-\alpha -2} (\xi\eta )^{-\alpha} e^{ -\frac{\xi^2+\eta^2}{2}} \big\{ ((\xi-\eta )^2-1)I_{\alpha}(\xi\eta )
  \\& \qquad\qquad\qquad\qquad\qquad
  -2 \xi\eta (I_{\alpha+1} (\xi\eta )-I_{\alpha}(\xi\eta ))-(2\alpha+1) \frac{\eta}{\xi } I_{\alpha+1}(\xi\eta  )\big\}
.\end{align*}
Combining these identities with \eqref{:A4} and \eqref{:A25}, and arguing as above, we obtain \eqref{:A3}.
For the case $k=2$, we also use  $\eta^2 /(1+\xi \eta) \le (1+|\xi-\eta|)^2$.

The estimates at $x=0$ follow by continuity, since the series expansion of $I_\alpha$ shows that $x\mapsto p_{t}^{0}(x,y)$ admits a smooth even extension to $\R$.
\end{proof}
\end{lemma}

Recall that $X$ is the $L_{\mathfrak{s},\alpha }^{(N)}$-diffusion. 
Let $\varphi (x)=\sqrt{2} \mathrm{arcsinh} \sqrt{x}$ and consider $Y_t:=\varphi(X_t)$.
Then, $Y:=\{Y_t\}_{t\ge 0}$ is the diffusion associated with the generator $L_0+ b(x) \partial _x$ 
on $(0,\infty)$, where $b(y)$ is the smooth bounded function on $[0,\infty)$ given by
\begin{align*}
b(y)=  \frac{(1-2N-\mathfrak{s})\sinh ^2(y /\sqrt{2}) +\alpha +\frac{1}{2} }{\sqrt{2} \sinh(y/\sqrt{2}) \cosh (y/\sqrt{2})}- \frac{\alpha+\frac{1}{2} }{y}
.\end{align*}

We use a parametrix construction adapted to the Bessel heat kernel; see, e.g., \cite{Fri64} for the standard parametrix method.
Set  $B:=B_x:=b(x) \partial_{x} $ and $H_t(x,y):=B_x p_{t}^{0}(x,y)$.
Define $  p_{t}^{(0)}(x,y):=  p_{ t}^{0}(x,y)$ and $  p_{t}^{(n+1)}(x,y) :=  (p^{(n)}*H)_{t}(x,y)$, where 
\begin{align*}
  (F*G)_t(x,y)=\int_0^t \int_0^\infty F_{t-s}(x,z) G_{s}(z,y) m_0(dz)ds
.\end{align*}

We can choose constants $\gamma, \Ct{c;A9}>0$ such that
\begin{align}\label{:A8}
p_{t}^{0}(x,y)
&\le \Ctref{c;A9}  p_{ \gamma t}^{0}(x,y), 
\\ \label{:A11}
|H_t(x,y)|
&\le \Ctref{c;A9} t^{-1/2} p_{ \gamma t}^{0}(x,y)
,\end{align}
by \eqref{:A2}, \eqref{:A3} and the boundedness of $b$. 
Hence, we have
\begin{align}\label{:A10}
|p_t^{(n)}(x,y)|
&\le  
 \frac{(\sqrt{\pi} \Ctref{c;A9}) ^{n+1}}{\Gamma(n/2+1)} p_{ \gamma t}^{0}(x,y) t^{n/2}
,\end{align}
for $n\ge 0$.
Indeed, this follows by induction in $n$, using \eqref{:A8}, \eqref{:A11}, and the Chapman-Kolmogorov equation for $p_{ \gamma t}^{0}(x,y)$.  
From \eqref{:A10}, the series
\begin{align*}
\tilde {p}_t(x,y)
:=
\sum_{n=0}^{\infty}p_t^{(n)}(x,y)
\end{align*}
converges absolutely and locally uniformly on $(0,T]\times [0,\infty)^2$.
Since each $p_t^{(n)}(x,y)$ is continuous on $(0,\infty) \ts [0,\infty)^2$ by induction, so is $\tilde {p}_t(x,y)$.


\begin{lemma}\label{l:A3}
$\tilde {p}_t(x,y)$ is the transition density of $Y$ with respect to $m_0$.
Furthermore, for $T>0$, there exist constants $\Ct{c;A13}, \Ct{c;A14}>0$ such that, for $0< t\le T $ and $x,y\ge 0$, we have 
\begin{align}\label{:A22}
\tilde{p}_{t}^{}(x,y) 
&
\le  \frac{\Ctref{c;A13}}{t^{1/2} (x+\sqrt{t})^{2\alpha+1}}\exp \Big(- \Ctref{c;A14} \frac{(x-y)^2}{t} \Big)
.\end{align}
\begin{proof}
Inequality \eqref{:A22} immediately follows from \eqref{:A2} and \eqref{:A10}.
Fix $f \in C_\infty( [0,\infty))$. 
Let $\{Q_t \}_{t\ge 0}$ be the Markov semigroup  for $Y$, and let $(P_{ t}^0 f)(x) =\int p_{t}^{0}(x,y) f(y )m_0(dy)$.
For $t>0$, $F_t (r) :=(P_t^0 f )(\sqrt{r} )$ is smooth up to $r=0$ by \eqref{:A1}.
To avoid the singular drift coefficient, consider the process $\{R_t\}_{t\ge 0}:=\{Y_t^2\}_{t\ge 0}$, which satisfies
\begin{align*}
  dR_t =2\sqrt{R_t}dB_t + 2\{\alpha+1+\sqrt{R_t} b(\sqrt{R_t})\} dt
.\end{align*}
Then, for fixed $t>0$ and $\ep >0 $, for $0\le s\le t-\ep$, It\^o's formula yields
\begin{align*}
dF_{t-s}(R_s)
&
=2 \sqrt{R_s}\partial_r F_{t-s} (R_s) dB_s +2 \sqrt{R_s} b(\sqrt{R_s}) \partial_r F_{t-s} (R_s) ds
\\&
=\partial_x (P_{t-s}^0f)(Y_s) dB_s+(B P_{t-s}^0f)(Y_s)ds 
.\end{align*}
Here we have used $\partial_t F_t=(2r \partial_r^2+2(\alpha+1)\partial_r)F_t$. 
Taking expectation, and using \eqref{:A3} to show the stochastic integral has zero expectation, we have 
\begin{align*}
Q_{t-\ep }(P_\ep^0 f)(x) =(P_t^0 f)(x)+\int_0^{t-\ep}Q_s(BP_{t-s}^0 f)(x )ds
,\end{align*}
for $x\ge 0$. 
Then, letting $\ep\to 0$, 
\begin{align}
  (Q_{t }f)(x) =(P_t^0 f)(x)+\int_0^{t}Q_{t-s}(H_s f)(x )ds
.\end{align}
Applying this repeatedly, we have 
\begin{align}\label{:A24}
Q_{t } f
=(P_t^0 f)+ \sum_{k=1}^{n-1 } (P^0 * H^{*k} )_t f +(Q * H^{*n} )_t f
.\end{align}
 By \eqref{:A11}, arguing as in the derivation of \eqref{:A10}, we have $\|  (Q * H^{*n} )_t f \|_\infty \to 0$ as $n\to\infty$. 
Hence \eqref{:A24} yields $Q_t f =\int_0^\infty \tilde{p}_t (x,y)f(y ) m_0(dy) $ by  the construction of $\tilde{p}_t$.
Thus, we complete the proof.
\end{proof}
\end{lemma}

Let $J(y):=\varphi(y)^{2\alpha+1 }\varphi' (y)$.
From \lref{l:A3} and $Y_t=\varphi (X_t)$, we have
\begin{align}\label{:43f}
 p_{\mathfrak{s},\alpha, t}^{(N)} (x,y)  =\tilde{p}_t(\varphi(x), \varphi(y)) J (y)  
.\end{align}

 \begin{lemma}\label{l:A2}
For every $t,y>0$,  $p_{\mathfrak{s}, \alpha, t}^{(N)}(\cdot,y) $ is smooth up to $x=0$. 
 Furthermore, for any $k\in \N$, $R>0$ and $t>0$, there exists a constant $\Ct{c;A10}:=\Ctref{c;A10}(k,t,R)$ such that, for $0<y\le R$, 
 \begin{align}\label{:A23}
 \sup_{x \ge 0}\Big| \frac{\partial^k }{\partial x^k } p_{\mathfrak{s}, \alpha, t}^{(N)}(x,y) \Big|
 &
 \le \Ctref{c;A10} J(y)
. \end{align} 
\begin{proof}
Let $\{Z_t \}_{t\ge 0}:=\{\log (1+X_t)\}_{t\ge 0}$ and define $\{ S_t\}_{t\ge 0}$ as the associated Markov semigroup.
The generator associated with $\{Z_t \}_{t\ge 0}$ is $(1-e^{-z})\partial_z^2+\{(\alpha+1)e^{-z}+(1-2N-\mathfrak{s})(1-e^{-z})\}\partial_z$, which satisfies the assumptions for \cite[Theorem 1.5]{Pop17}.
Let $h_y(z):= p_{\mathfrak{s}, \alpha, t/2}^{(N)}(e^z-1,y) =\tilde{p}_{t/2}(\varphi(e^z-1),\varphi(y)) J(y)$.
Then, we have $h_y \in C_\infty ([0,\infty) )$ and $\sup_{0< y\le R} \| h_y \|_\infty /J(y) <\infty$ from \lref{l:A3}.
 Hence, applying \cite[Theorem 1.5]{Pop17} to $ S_{t/2} h_y  (z )$, we have 
 \begin{align*}
  \max_{0\le j\le k} \sup_{z\ge 0} |\partial_z^j S_{t/2} h_y  (z ) |\le \mathrm{Const}_{k,t}\| h_y\|_\infty \le  \Ctref{c;A10} J(y)
. \end{align*}
Thus, using $ p_{\mathfrak{s}, \alpha, t}^{(N)}(x,y)=( S_{t/2} h_y ) ( \log(1+x) )$, we complete the proof.
\end{proof}
\end{lemma}

\section*{Acknowledgement}

We are deeply grateful to Grigori Olshanski for his valuable suggestion that led to the result in \sref{s:6}.
A.B. is a winner of the “BASIS” Foundation Competition in Mathematics and theoretical physics and is deeply grateful to the Jury and the sponsors.

\section*{Funding}
A.B.’s  research was  supported by the Ministry of Science and Higher Education of the Russian Federation, project NO. 075-15-2024-631.
K.Y. is supported by JSPS KAKENHI Grant Numbers JP21K13812 and JP25K17268.

\bibliographystyle{acm} 
\bibliography{reference}

\end{document}